\documentclass[onefignum,onetabnum]{siamonline220329}

\usepackage{lipsum}
\usepackage{amsfonts}
\usepackage{graphicx}
\usepackage{epstopdf}
\ifpdf
  \DeclareGraphicsExtensions{.eps,.pdf,.png,.jpg}
\else
  \DeclareGraphicsExtensions{.eps}
\fi

\usepackage{enumitem}
\setlist[enumerate]{leftmargin=.5in}
\setlist[itemize]{leftmargin=.5in}

\newsiamremark{remark}{Remark}
\newsiamremark{hypothesis}{Hypothesis}
\crefname{hypothesis}{Hypothesis}{Hypotheses}
\newsiamthm{claim}{Claim}

\headers{Go-UQ for Inverse Problems via VED Networks}{B. Afkham, J. Chung, and M. Chung}

\DeclareMathOperator{\diag}{diag}

\ifpdf
\hypersetup{
  pdftitle={Goal-oriented Uncertainty Quantification for Inverse Problems via Variational Encoder-Decoder Networks},
  pdfauthor={B. Afkham, J. Chung, and M. Chung}
}
\fi

\usepackage{tiaArticleStyle}

\usepackage{makecell}
\usepackage{subcaption}
\usepackage{soul}

\graphicspath{{media/}}
\bibliography{references.bib}

\newcommand{\jmc}[1] {{\color{blue} #1}}
\newcommand{\tia}[1] {{\color{blue} #1}}

\newtheorem{assumption}{Assumption}

\begin{document}

\title{Goal-oriented Uncertainty Quantification for Inverse Problems via Variational Encoder-Decoder Networks\thanks{Submitted to the editors \today.
\funding{This work was partially supported by the National Science Foundation (NSF) under grants DMS-1654175 and 2026841 (J. Chung), DMS-2152661(M. Chung) and partially supported by The Villum Foundation Grant No. 25893 (B M Afkham) }}}

\author{Babak Maboudi Afkham\thanks{DTU Compute, Department of Applied Mathematics and Computer Science, Technical University of Denmark, Lyngby, Denmark}
\and Julianne Chung\thanks{Department of Mathematics, Emory University, Atlanta, GA 30322, USA  (\email{jmchung@emory.edu}, \url{http://www.math.emory.edu/\string~jmchung/})}
\and Matthias Chung\thanks{Department of Mathematics, Emory University, Atlanta, GA 30322, USA}}

\maketitle

\begin{abstract}
In this work, we describe a new approach that uses variational encoder-decoder (VED) networks for efficient
goal-oriented uncertainty quantification for inverse problems.
Contrary to standard inverse problems, these approaches are \emph{goal-oriented} in that the goal is to estimate some quantities of interest (QoI) that are functions of the solution of an inverse problem, rather than the solution itself. Moreover, we are interested in computing uncertainty metrics associated with the QoI, thus utilizing a Bayesian approach for inverse problems that incorporates the prediction operator and techniques for exploring the posterior.  This may be particularly challenging, especially for nonlinear, possibly unknown, operators and nonstandard prior assumptions. We harness recent advances in machine learning, i.e., VED networks, to describe a data-driven approach to large-scale inverse problems. This enables a real-time goal-oriented uncertainty quantification for the QoI. One of the advantages of our approach is that we avoid the need to solve challenging inversion problems by training a network to approximate the mapping from observations to QoI.  Another main benefit is that we enable uncertainty quantification for the QoI by leveraging probability distributions in the latent space.  This allows us to efficiently generate QoI samples and circumvent complicated or even unknown forward models and prediction operators. Numerical results from medical tomography reconstruction and nonlinear hydraulic tomography demonstrate the potential and broad applicability of the approach.
\end{abstract}



\section{Introduction} \label{sec:introduction}
Inverse problems arise in many scientific applications, and there has been a significant amount of research on the development of theory and computational methods for solving inverse problems.  However, solving an inverse problem is often just one step in a multi-step process, where the solution of the inverse problem is a tool that is used for some end goal (e.g., predictions that are used for optimal design or control). Standard approaches for goal-oriented inverse problems first seek a solution to the inverse problem and, informed by this solution, subsequently compute quantities of interest (QoI).  While in some cases, estimation of these QoI is the end goal, there are many scenarios where uncertainties for the QoI provide critical information, e.g., for further decision-making. In this work, we are interested in \emph{goal-oriented uncertainty quantification} (go-UQ), which combines goal-oriented inverse problems with uncertainty quantification for the QoI.  In particular, we provide an efficient end-to-end computational approach that can be used to directly estimate the QoI and their uncertainties.

\medskip

Consider a \emph{goal-oriented inverse problem}
of the form
\begin{equation} \label{eq:detmodel}
    \begin{aligned}
    & \bfb = F(\bfy) + \bfe, \\
    & \bfx = G(\bfy),
    \end{aligned}
\end{equation}
where $\bfb \in \bbR^m$ contains the observations, $\bfy\in \bbR^n$ contains unknown true parameters (e.g., the true image or state parameters), $F:\bbR^n \to \bbR^m$ is the parameter to observable map, and $\bfe\in \bbR^m$ represents additive noise. In classical inverse problems, the goal is to compute approximations of the unknowns in $\bfy$, given observations $\bfb$ and knowledge of the forward operator $F$, \cite{hansen2010discrete}.   Here, we are interested in \emph{goal-oriented} inverse problems \cref{eq:detmodel} where we aim to estimate some QoI in $\bfx\in \bbR^q$ that may be determined by a \emph{prediction operator} $G:\bbR^n \to \bbR^q$, mapping the unknown $\bfy$ onto $\bfx$.

Some motivating examples of goal-oriented inverse problems include carbon sequestration \cite{lieberman2013goal}, modeling transport of a contaminant in an urban environment \cite{attia2018goal}, and detecting the locations of discontinuities in reconstructed CT images \cite{dahl2017computing,afkham2023uncertainty}. In many applications, the prediction operator may be simply taking a maximum, minimum, average, or some combination of the unknown parameters in $\bfy$, \cite{lieberman2013goal,devore2019computing,spantini2017goal, attia2018goal}. However, QoIs may also be related to $\bfy$ by more complex mappings where the prediction operator may not be well-defined or may rely on expert opinion (e.g., manual classification or segmentation of images). Regardless, the QoI typically represents some low dimensional property of the unknown $\bfy$, i.e., $q \ll n$.

Goal-oriented inverse problems suffer from the same obstacles as classical inverse problems such as ill-posedness and large parameter dimensions.
Further challenges arise when considering UQ for the QoI due to the additional prediction operator $G$.
Recent works on goal-oriented inverse problems consider modifications of the inversion process or the data acquisition process to take into account the output QoI.  For example, in goal-oriented inference, end goals (i.e., output QoI) are incorporated into the inference process thereby leading to efficient inference-for-prediction algorithms \cite{lieberman2013goal,lieberman2014nonlinear}. These approaches rely on estimating the uncertainty of some predictive QoI based on the solution of an inverse problem. However, the main difference is that these approaches target the most relevant parameters for prediction (i.e., modifying the inference process to account for output quantities of interest), and we are interested in estimating uncertainties of the QoI as the end goal.  Goal-oriented approaches have also been considered in the context of optimal approximations \cite{spantini2017goal} and optimal design of experiments \cite{attia2018goal}.  These approaches use the output QoI, which are often linear mappings of the unknowns, to inform dimensionality reduction techniques and experimental design (e.g., optimal placement of sensors) respectively.

\paragraph{Overview of contributions} In this work, we provide a novel computational framework for directly estimating the QoI and their uncertainties, given observations.  We assume that sample pairs $\left\{(\bfb^j, \bfx^j) \right\}_{j=1}^J$ of observations and QoI are readily available and use a supervised learning approach to train a \emph{variational encoder-decoder} (VED) network to approximate the mapping from observations $\bfb$ to QoI $\bfx$. This approach confers several advantages:
\begin{itemize}
    \item We bypass the computationally challenging inversion process of estimating $\bfy$ by establishing an end-to-end framework that maps the observations $\bfb$ directly to the QoI $\bfx$. Although one could consider a two-step approach where a network is used to approximate the inversion process from $\bfb$ to $\bfy$ \cite{li2020nett}, followed by the QoI computation, there are significant challenges due to the potentially large number of unknowns in $\bfy$. Dimension reduction techniques have been considered, e.g., to train a deep neural network (DNN) to approximate the mapping from observation data to desired (e.g., regularization) parameters \cite{afkham2021learning,liu2021machine}; however, these approaches do not allow for UQ of the QoI.
    \item By leveraging VED networks, we can perform efficient UQ of the QoI.  More specifically, we use a VED network to create a mapping from the observation to some latent space, enforce a predetermined distribution on the unobserved latent variables, and then create a mapping from the latent space to the target data distribution (e.g., the QoI).
    The latter half of the VED defines a generating function that maps independent and identically distributed (i.i.d.) samples from the latent space variables to the QoI.  Once the network is trained, UQ for the QoI can be performed very efficiently in an online phase, via forward propagation of the observation $\bfb$ through the VED, where the output of the VED is a set of samples that can be used for the statistical description of the QoI $\bfx$.
    \item One key advantage of our approach is that it is data-driven, in that training data are used to circumvent complicated or even unknown forward models and prediction operators. Most goal-oriented approaches use linear prediction operators $G$. Our approach is not restricted to and does not require knowledge of $G$ (or the forward operator $F$).
\end{itemize}

An overview of the paper is as follows.  In \Cref{sec:background} we provide a brief description of the Bayesian formulation for goal-oriented inverse problems and describe some basics of VED networks for readers unfamiliar with these machine learning techniques.  Then in  \Cref{sec:ved_goUQ}, we describe various details about using VED networks for go-UQ for inverse problems. Numerical experiments provided in \Cref{sec:numerics} illustrate the benefits of our approach for various example applications and goal-oriented problem setups.  Conclusions are provided in \Cref{sec:conclusions}.

\section{Background}
\label{sec:background}
Recall that for goal-oriented inverse problems \cref{eq:detmodel}, the standard approach is to first estimate the unknowns in $\bfy$, given $\bfb$ and $F$, and second to compute the QoI $\bfx$ from $\bfy$.  However, standard approaches often do not reveal the uncertainties in this estimation. A popular approach that enables quantification of uncertainties (on $\bfy$, $\bfx$, or both) is to provide a probabilistic formulation for \cref{eq:detmodel}.
That is, in order to provide information about the accuracy or uncertainties of the estimate, we are motivated to consider inverse problems in the framework of Bayesian statistics. In \Cref{sub:Bayesian}, we review the Bayesian formulation of an inverse problem with particular emphasis on prediction for the goal-oriented problem.
Then, in \Cref{sub:NN} we describe tools from machine learning and their use within the inverse problems community.  In particular, we describe VED networks that will be used in \Cref{sec:ved_goUQ} for go-UQ.

\subsection{Bayesian inverse problem}
\label{sub:Bayesian}
In this section, we describe a statistical approach to solving inverse problems. This approach enables us to include uncertainties, e.g., in measurement, into our model. In addition, the statistical approach provides a natural method for incorporating our prior knowledge of the parameters. Contrary to deterministic inverse problems where the goal is to compute a point estimate, the solution of a Bayesian inverse problem is a random variable and thus has a distribution \cite{calvetti2007introduction}.
More specifically, in the statistical framework, all unknown values are modeled as random variables.  Let $B, X, Y,$ and $E$ be random vectors representing the observational data, the QoI, the unknown state, and additive noise, respectively. Then, the Bayesian interpretation of \cref{eq:detmodel} is given by
\begin{align}
    & B = F(Y) + E, \label{eq:stochmodel1}\\
    & X = G(Y).     \label{eq:stochmodel2}
\end{align}
Note that since $Y$ and $E$ are unknown, we represent their uncertainties with a probability distribution constructed from prior knowledge. In the following, assuming for simplicity of illustration that all densities exist, the solution to the inverse problem \cref{eq:stochmodel1} is given by the posterior density function, conditioned on the observations $\bfb$, i.e.,
 \begin{equation} \label{posterior}
     \pi_{\rm post}(\bfy \mid \bfb) \propto \pi_{\rm like}(\bfb \mid \bfy)\, \pi_{\rm prior}(\bfy).
 \end{equation}
 Since we are interested only in the distribution of the QoI $X$ rather than those of the state parameters $Y$ themselves, the posterior distribution is propagated through to the QoI $X$, i.e.,
 \begin{equation}
 \label{eq:posteriorpredictive}
    X \mid B = G(Y \mid B).
 \end{equation}
In the context of goal-oriented inference problems, this is often referred to as the \emph{posterior predictive}, \cite{lieberman2013goal}. The prediction is a random variable and, hence, is characterized by a distribution. We denote its density by $\pi_{\rm pred}(\bfx \mid \bfb)$.  Given a prior and some observation, samples from the posterior $Y \mid B$ can be passed through the prediction process, and the resulting samples can be used to describe the posterior predictive and provide uncertainties for the QoI.\\

In the special case where the parameter-to-observable map is linear $F(\bfy) = \bfA \bfy$, the noise $E$ and the prior $Y$ are Gaussian, i.e., $E\sim\calN(\bfzero, \bfGamma_{\rm noise})$ and $Y\sim\calN(\widebar\bfy, \bfGamma_{\bfy})$ with $\bfGamma_{\rm noise}$ and $\bfGamma_{\bfy}$ symmetric and positive definite covariance matrices, the posterior \cref{posterior} is also a Gaussian and is given by $\calN(\bfy_{\rm post}, \bfGamma_{\rm post})$ where
\begin{equation}
    \bfGamma_{\rm post} = \left(\bfA\t \bfGamma_{\rm noise}^{-1} \bfA + \bfGamma_\bfy^{-1}\right)^{-1}
    \quad \mbox{and}\quad
     \bfy_{\rm post} = \bfGamma_{\rm post}\left(\bfGamma_\bfy^{-1}\widebar\bfy + \bfA\t \bfGamma_{\rm noise}^{-1}\bfb\right).
\end{equation}
If, additionally, the prediction operator is linear, i.e., $G(\bfy) = \bfP \bfy$ for appropriately sized matrix $\bfP$, then the posterior predictive \cref{eq:posteriorpredictive} is also a Gaussian and is given by $\calN(\bfx_{\rm pred}, \bfSigma_{\rm pred})$ where
\begin{equation}
    \bfx_{\rm pred} = \bfP  \bfy_{\rm post}
\quad \mbox{and}\quad \bfSigma_{\rm pred} = \bfP \bfGamma_{\rm post} \bfP\t.
    \end{equation}
In the linear setting, a computationally efficient inference-for-prediction method was described in \cite{lieberman2013goal} that exploits underlying balanced truncation model reduction for dimension reduction.  The posterior predictive was also used for the goal-oriented optimal design of experiments in \cite{attia2018goal}.

So far, we have described goal-oriented approaches, and for the special case described above, efficient tools have been developed for exploring the posterior predictive.
In this work, we exploit VED networks for go-UQ, resulting in a data-driven approach that can exploit recent developments in machine learning, is generalizable, and works for complex and even unknown prediction operators.

\subsection{Learning from data} \label{sub:NN}
Machine learning, in particular, deep neural networks has become an important tool in the inverse problems community.  In this section, we provide some background on neural networks and describe the training process, given some training data.
We describe how these networks work by  starting with a description of DNNs, then describing encoder-decoder networks.

Let us begin with the assumption that there exists some target function $\widebar\bfPhi:\bbR^{m} \to \bbR^q$ that maps input vector $\bfb \in \bbR^{m}$ to a target output vector $\bfx \in \bbR^{q}$.

\paragraph{Deep neural networks (DNNs)}
\emph{Deep feed-forward networks} or \emph{deep neural networks} (DNNs) are among the many classic machine learning modeling tools. These networks are comprised of a composition of multiple layers of computation. In a feed-forward network, the layers together with their connections are typically represented by a directed graph, indicating the flow of information from the input to the target.

For simplicity, consider a feed-forward DNN with $K$ layers for an input vector $\bfb \in \bbR^{m}$ and a target output vector $\bfx \in \bbR^{q}$ as
\begin{equation}
     \bfPhi(\bfb; \bftheta) = \bfvarphi_{K+1}(\bftheta_{K+1})\circ \dots \circ\bfvarphi_{1}(\bftheta_{1})(\bfb) = \bfx,
\end{equation}
where $\bfx$ is the output of the network, $\circ$ denotes the component-wise composition of functions $\bfvarphi_k:\bbR^{r_{k-1}} \times \bbR^{p_k} \to \bbR^{r_k}$ for $k = 1, \ldots, K+1$  ($r_0 = m$, and $r_{K+1} = q$). The vector $\bftheta = \begin{bmatrix} \bftheta_1\t, & \ldots, & \bftheta_{K+1}\t\end{bmatrix}\t \in\bbR^p$ is a composition of layer specific parameters $\bftheta_k$ defining the so-called \emph{weights} $\bfW_k$ and \emph{biases} $\bfd_k$, i.e., $\bftheta_k = \begin{bmatrix} \vec{\bfW_k}\t & \bfd_k\t \end{bmatrix}\t$, where $\bfW_k \in\bbR^{r_
k \times r_{k-1}}$ and $\bfd_k \in \bbR^{r_k} $. The functions $\bfvarphi_k$ are given by
\begin{equation}
    \bfb_k = \bfvarphi_k(\bftheta_k)(\bfb_{k-1}) = \bfs_k( \bfW_k \bfb_{k-1} + \bfd_k ),
\end{equation}
with $\bfb_0 = \bfb$ and
where $\bfs_k:\bbR^{r_k} \to \bbR^{r_k}$ are typically nonlinear \emph{activation functions}, mapping input arguments point-wise onto outputs with limiting range. Note that for the \emph{output layer} $\bfvarphi_{K+1}(\bftheta_{K+1})$, we assume a linear transformation with no bias term, $\bfvarphi_{K+1}(\bftheta_{K+1})(\bfb_K) = \bfW_{K+1} \bfb_{K}$.\\

Given input-target pairs $\left\{(\bfb^j,\bfx^j)\right\}_{j = 1}^J$, the aim of supervised learning is to train the neural network, i.e., computing network parameters $\bftheta$, such that  $\bfPhi(\bfb;\bftheta)\approx \widebar\bfPhi(\bfb)$. For this, a \emph{loss function} determines the performance of a network by comparing $\bfPhi(\bfb^j;\bftheta)$ and $\bfx^j$ for a particular $\bftheta$.
An optimization method is used to find appropriate network parameters $\bftheta$ by minimizing a suitable loss function.

DNNs have been used for solving inverse problems (e.g., for regularization parameter selection \cite{afkham2021learning,liu2021machine} or for learning regularization functionals \cite{de2016machine,hammernik2018learning}). Furthermore, DNNs have been adopted for learning the inversion process, i.e., the mapping from the data space $\bfb$ to the inverse solution $\bfy$, \cite{li2020nett}. Such approaches could be extended for the prediction or estimation of the QoI.

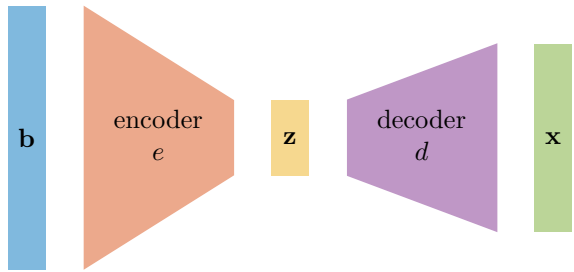
\begin{figure}[bthp]
    \centering
    \begin{tikzpicture}
     
	\node[fill=matlab1!50, minimum width=0.5cm, minimum height=3.5cm] (b) at (0,0) {$\bfb$};
	
	\draw[fill=matlab2!50,draw=none] ([xshift=0.5cm]b.north east) -- ([xshift=2.5cm,yshift=0.5cm]b.east) -- ([xshift=2.5cm,yshift=-0.5cm]b.east) -- ([xshift=0.5cm]b.south east) -- cycle; 
	\node at (1.75,0) {$\begin{matrix}\mbox{encoder} \\ e \end{matrix}$};
	
	\node[fill=matlab3!50, minimum width=0.5cm, minimum height=1.0cm] (Z) at (3.5cm,0) {$\bfz$};
	
	\draw[fill=matlab4!50,draw=none] ([xshift=0.5cm]Z.north east) -- ([xshift=2.5cm,yshift=0.75cm]Z.north east) -- ([xshift=2.5cm,yshift=-0.75cm]Z.south east) -- ([xshift=0.5cm]Z.south east) -- cycle;
	\node at (5.25,0) {$\begin{matrix}\mbox{decoder} \\ d \end{matrix}$};
	
	\node[fill=matlab5!50, minimum width=0.5cm, minimum height=2.5cm] (Xp) at (7,0) {$\bfx$};

\end{tikzpicture}
    \caption{Encoder-decoder network mapping input vector $\bfb$ to output vector $\bfx$.  The encoder maps the input vector $\bfb$ to the latent variables $\bfz$ and the decoder maps the latent variables to the output $\bfx$.}
    \label{fig:EDnetwork}
\end{figure}

\paragraph{Encoder-Decoder Networks}
DNNs can be implemented within an Encoder-Decoder (ED) architecture. Such networks consist of the concatenation of two networks (the encoder and the decoder). The encoder network is comprised of layers that successively reduce the dimension of the input vector. The output of the encoder is referred to as the \emph{latent} vector $\bfz \in \bbR^\ell$ that describes a representation of the input. Then, the decoder network increases the dimension through the network to match the target dimension.  More specifically, the two main parts can be defined using transformations:
\begin{itemize}
    \item the encoder network $e:\bbR^m \times \bbR^{p^{\rm e}}\to \bbR^\ell$ with network parameters $\bftheta^{\rm e} \in \bbR^{p^{\rm e}}$ maps an input $\bfb$ to the latent variables $\bfz$,
    \begin{equation}
        \bfz = e(\bfb;\,\bftheta^{\rm e}),
    \end{equation}
    \item the decoder network $d:\bbR^\ell\times \bbR^{p^{\rm d}} \to \bbR^q$ with network parameters $\bftheta^{\rm d} \in \bbR^{p^{\rm d}}$ maps the latent variables to the output $\bfx$,
    \begin{equation}
        \bfx = d(\bfz;\, \bftheta^{\rm d}).
    \end{equation}
\end{itemize}
The entire ED network then takes the form $\bfPhi_{\rm ed}(\bfb; \, \bftheta^{\rm e}, \bftheta^{\rm d}) := d(e(\bfb; \, \bftheta^{\rm e})\, ;\, \bftheta^{\rm d})$. See \cref{fig:EDnetwork} for a visual representation of an ED network.

Analogous to DNNs, the aim is to find parameters $\bftheta^{\rm e}$ and $\bftheta^{\rm d}$ that minimize some \emph{loss function}. A common approach is to seek network parameters $\bftheta^{\rm e}, \bftheta^{\rm d}$ to minimize the mean squared error, i.e.,
\begin{equation}
\min_{\bftheta^{\rm e}, \bftheta^{\rm d}} \ \sum_{j = 1}^J \norm[2]{\bfx^j - \bfPhi_{\rm ed}(\bfb^j; \, \bftheta^{\rm e}, \bftheta^{\rm d})}^2.
\end{equation}

A commonly-used ED network is the \emph{autoencoder}, which aims to map an input $\bfb$ onto itself, i.e., $\bfPhi_{\rm ed}(\bfb; \, \bftheta^{\rm e}, \bftheta^{\rm d}) \approx \bfb$. The encoder represents the input in a lower dimensional space while the decoder attempts to reconstruct the input from the low dimensional representation $\bfz$. The flow of data through a low-dimensional latent space forces the autoencoder to find a low-dimensional representation for a set of data, by which insignificant data is eliminated. Hence, such networks are classically used in dimensionality and noise reduction, and data compression applications, see e.g., \cite{hinton2006reducing,salakhutdinov2007restricted,torralba2008small,wang2016auto}.

\paragraph{Variational Encoder-Decoder (VED) Networks}
Standard ED networks do not have \emph{generative} capability, meaning a random sample from the latent space $\bfz^j \in \bbR^\ell$ may not necessarily generate a meaningful target sample $d(\bfz^j;\bftheta^{\rm d})$. Variational Encoder-Decoder (VED) networks address this issue by enforcing a predetermined structure/distribution on the latent variables.
This is particularly relevant in the context of variational autoencoders (VAEs), where the goal is to have the output match the input.  Here, by sampling from the predetermined latent distribution, the decoder $d$ provides a generating function that can be used to simulate or generate new data with desired features \cite{kingma2013auto}.

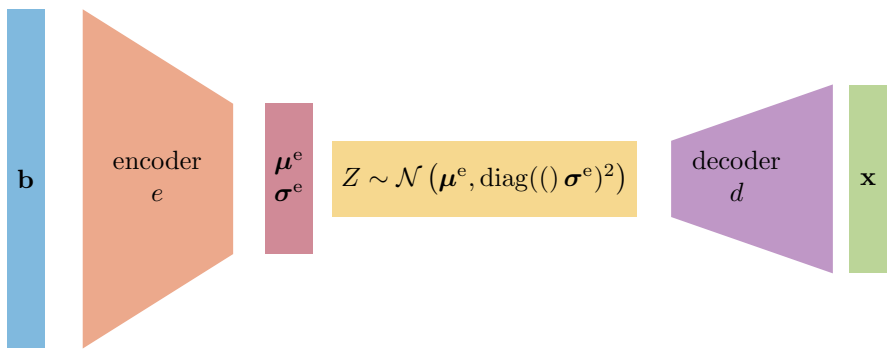
\begin{figure}[bthp]
    \centering
    \begin{tikzpicture}
     
	\node[fill=matlab1!50, minimum width=0.5cm, minimum height=4.5cm] (X) at (0,0) {$\bfb$};
	
	\draw[fill=matlab2!50,draw=none] ([xshift=0.5cm]X.north east) -- ([xshift=2.5cm,yshift=1.0cm]X.east) -- ([xshift=2.5cm,yshift=-1.0cm]X.east) -- ([xshift=0.5cm]X.south east) -- cycle; 
	\node at (1.75,0) {$\begin{matrix}\mbox{encoder} \\ e \end{matrix}$};
	
	\node[fill=matlab7!50, minimum width=0.5cm, minimum height=2.0cm] (sa) at (3.5cm,0) {$\begin{matrix}\bfmu^{\rm e} \\ \bfsigma^{\rm e}\end{matrix}$};
	
	\node[fill=matlab3!50, minimum width=0.5cm, minimum height=1.0cm] (mu) at (6.1cm,0) {$Z \sim \calN\left(\bfmu^{\rm e},\diag(\bfsigma^{\rm e})^2\right)$};

	\draw[fill=matlab4!50,draw=none] ([xshift=0.45cm]mu.north east) -- ([xshift=2.6cm,yshift=0.75cm]mu.north east) -- ([xshift=2.6cm,yshift=-0.75cm]mu.south east) -- ([xshift=0.45cm]mu.south east) -- cycle;
	\node at (9.45,0) {$\begin{matrix}\mbox{decoder} \\ d \end{matrix}$};
	
	\node[fill=matlab5!50, minimum width=0.5cm, minimum height=2.5cm] (Xp) at (11.2,0) {$\bfx$};
     
\end{tikzpicture}
    \caption{Variational Encoder-Decoder network mapping input vector $\bfb$ to target vector $\bfx$. The encoder maps the input vector $\bfb$ to the latent distribution defined by mean $\bfmu^{\rm e}$ and standard deviations $\bfsigma^{\rm e}$. The latent variables $\bfz$ are obtained by sampling from the latent distribution, and the decoder maps the latent variables to the target $\bfx$.}
    \label{fig:VEDnetwork}
\end{figure}

In a VED network, the input $\bfb$ is passed through an encoder  $e(\bfb;\,\bftheta^{\rm e})$ that depends on parameters $\bftheta^{\rm e}$. This is similar to the ED networks where we reduce the dimensions of the input. However,
the output of the encoder is a mean vector $\bfmu^{\rm e} \in \bbR^\ell$ and a vector of standard deviations $\bfsigma^{\rm e} \in \bbR^\ell$. For simplicity, we assume that $\bfmu^{\rm e}$ and $\bfsigma^{\rm e}$ define a Gaussian distribution for the latent vector $\bfz$, i.e.,

\begin{equation}
    Z \sim \calN\left(\bfmu^{\rm e},\diag{(\bfsigma^{\rm e})}^2 \right).
\end{equation}
Then a sample $\bfz$ is fed to the decoder $d(\bfz;\,\bftheta^{\rm d})$ that depends on parameters $\bftheta^{\rm d}$.  Notice that contrary to the ED network, the latent vector $\bfz$ is not the output of the encoder but rather a sample from the latent distribution. See \cref{fig:VEDnetwork} for a visual representation of a VED network.

Training a VED network requires determining network parameters $\bftheta^{\rm e}$ and $\bftheta^{\rm d}$ that minimize some objective or loss.
However, contrary to DNNs and EDs, for VEDs there are two main components to consider.  First, we wish to reconstruct an output matching a target by minimizing a training loss, e.g., for a sample $\{\bfx^j,\bfb^j\}$
\begin{equation}
    \calD_{\rm train}^j(\bftheta^{\rm e},\bftheta^{\rm d}) =  \norm[2]{\bfx^j - d(\bfz;\ \bftheta^{\rm d})}^2
    \end{equation}
where
\begin{equation}
Z \sim \calN\left(\bfmu^{\rm e}, \diag{(\bfsigma^{\rm e})}^2\right) \quad \mbox{with} \quad \begin{bmatrix} \bfmu^{\rm e} \\ \bfsigma^{\rm e} \end{bmatrix} =  e(\bfb^j; \, \bftheta^{\rm e}).
\end{equation}
Second, we aim to drive the latent space distribution toward the isotropic Gaussian (standard Gaussian with zero mean and unit variance), which can be obtained by minimizing a \emph{similarity} loss between distributions.
For example, we can define the similarity loss for the VED network as
\begin{equation}
    \calD_{\rm sim}^j(\bftheta^{\rm e}) = \calD_{\rm kl}  \big( \calN(\bfmu^{\rm e}, \diag{(\bfsigma^{\rm e})}^2)\,  \mid\mid \, \calN(\bfzero, \bfI_\ell) \big),
\end{equation}
where $\calD_{\rm kl}$ is the Kullback–Leibler (KL) divergence or KL loss \cite{kullback1951information}, and the objective is to find encoder and decoder parameters $\bftheta^{\rm e}$ and $\bftheta^{\rm d}$ that minimize the total loss \cite{kingma2013auto},
\begin{equation}
\label{eq:ELBO}
    \sum_{j = 1}^J \calL^j(\bftheta^{\rm e}, \bftheta^{\rm d})
    \quad \mbox{where} \quad
    \calL^j(\bftheta^{\rm e}, \bftheta^{\rm d}) = \calD_{\rm train}^j(\bftheta^{\rm e},\bftheta^{\rm d}) + \calD_{\rm sim}^j(\bftheta^{\rm e}).
\end{equation}

\emph{Variational autoencoders} are special cases of VED networks, and they have been used in various contexts for solving inverse problems, especially in imaging \cite{lucas2018using, bora2017compressed,peng2020solving}.
These variational autoencoders were also used in the context of uncertainty quantification in \cite{goh2019solving}, where they are used to learn the parameter to observable map as well as the inverse problem solver.
However, since the input and output dimensions must be the same for variational autoencoders and adversarial autoencoders, these approaches are too restrictive to be considered for goal-oriented UQ.

\section{Variational Encoder-Decoder (VED) networks for goal-oriented UQ}
\label{sec:ved_goUQ}
In this section, we describe a VED network that can be used to approximate the posterior predictive \cref{eq:posteriorpredictive}. We establish a relationship between random variables $B$ and $X$. The maps in this graph are approximated using a VED network.
In \Cref{sub:learnIPT}, we argue how the relation between $B$ and $X$ results in a loss function similar to the one expressed in \cref{eq:ELBO}. We describe how training data can be used to estimate network parameters.
Once the VED is trained and we have a new observation, we describe in \Cref{sub:UQforQoI} a sampling approach to perform UQ for the QoI.

Recall that the density function for the QoI $\pi(\bfx)$ is, in general, unknown. We aim to approximate this distribution, as well as the density of the posterior predictive $\pi_{\rm pred}(\bfx \mid \bfb)$.
In this approach we assume there exists a latent variable with a convenient density function, e.g., a Gaussian density function, and consider a
map $\widebar d:\bbR^\ell \to \bbR^q$ such that $\bfx = \widebar d(\bfz)$.
The benefit of such a map \cite{lemaire2013structural} is that it reduces the task of exploring an unknown and complex density function, e.g., $\pi(\bfx)$, to exploring a computationally convenient density function. However, the complexity of $\pi(\bfx)$ is still encoded in the nonlinear map $\widebar d$.

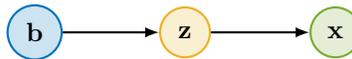
\begin{figure}[ht]
    \centering
\begin{center}
    \resizebox{0.3\textwidth}{!}{
\begin{tikzpicture}
[block/.style ={circle, draw=matlab1, thick, fill=matlab1!20, text width=2ex,align=center, minimum height=1em}, scale = 1]
\node[block,draw=matlab1,fill=matlab1!20] at (0, 0) (b) {$\bfb$};
\node[block,draw=matlab3,fill=matlab3!20] at (2, 0) (z) {$\bfz$};
\node[block,draw=matlab5,fill=matlab5!20] at (4, 0) (x) {$\bfx$};
\draw[-latex, black, thick] (b) -- (z);
\draw[-latex, black, thick] (z) -- (x);
\end{tikzpicture}
}
\end{center}
    \caption{Graphical model representing the connection between the observation $\bfb$, the latent variables $\bfz$ and the QoI $\bfx$.}
    \label{fig:graph}
\end{figure}

We now relate $\pi(\bfz)$ to the observation density $\pi(\bfb)$ by introducing another map $\widebar e:\bbR^m \to \bbR^\ell$, such that $\bfz = \widebar e(\bfb)$. Note that, in this case, the density function $\pi(\bfb)$ is unknown. The composite map $\widebar d \circ \widebar e$ establishes a connection between $B$ and $X$. We represent the connection between the random variables $B$, $Z$, and $X$ in the graphical model in \cref{fig:graph}.

\begin{assumption} \label{assump:bias}
We assume that $Z$ parameterizes $X$ such that $X$ and $Z$ are indistinguishable to $B$, i.e.,
\begin{equation} \label{eq:bias}
    \pi(\bfb\mid \bfx,\bfz) = \pi(\bfb\mid\bfz) = \pi(\bfb\mid\bfx).
\end{equation}
\end{assumption}
In other words, \Cref{assump:bias} suggests that the parameterization of the prior distribution does not affect the data generation process. Now we check if the latent variable $Z$ introduces a bias in the posterior predictive. We first construct an expression for the joint distribution $\pi(\bfb,\bfx,\bfz)$ following \Cref{assump:bias}. We have
\begin{equation}
        \pi(\bfb\mid\bfx,\bfz) = \frac{\pi(\bfb,\bfx,\bfz)}{\pi(\bfx,\bfz)} \qquad \mbox{and} \qquad \pi(\bfb\mid\bfx) = \frac{\pi(\bfb,\bfx)}{\pi(\bfx)}.
\end{equation}
Combining this with \cref{eq:bias} yields
\begin{equation} \label{eq:joint}
    \pi(\bfb,\bfx,\bfz) = \frac{\pi(\bfb,\bfx)\pi(\bfx,\bfz)}{\pi(\bfx)}.
\end{equation}
Now we evaluate the bias in the posterior predictive in the presence of the latent variable $\bfz$
\begin{equation}
    \pi(\bfx\mid\bfb,\bfz) = \frac{\pi(\bfx,\bfb,\bfz)}{\pi(\bfb,\bfz)} = \frac{\pi(\bfb,\bfx)\pi(\bfx,\bfz)}{\pi(\bfx)\pi(\bfb,\bfz)},
\end{equation}
where we used \cref{eq:joint} in the last equality. Multiplying both the numerator and the denominator with $\pi(\bfb)$ yields
\begin{equation}
    \pi(\bfx\mid\bfb,\bfz) =  \frac{\pi{(\bfb)}}{\pi(\bfb,\bfz)} \frac{\pi(\bfx,\bfz)}{\pi(\bfx)} \frac{\pi(\bfb,\bfx)}{\pi(\bfb)} = \frac{\pi(\bfz\mid\bfx)}{\pi(\bfz\mid\bfb)} \pi_{\rm pred}(\bfx\mid\bfb).
\end{equation}
Here we used the definition of a conditional distribution. Note that the right-hand-side contains the posterior predictive $\pi_{\rm pred}(\bfx\mid\bfb)$, and
 the ratio $\pi(\bfz\mid\bfx)/\pi(\bfz\mid\bfb)$ does not equal to one for general maps $\widebar d$ and $\widebar e$. However, we require $\widebar d$ and $\widebar e$ to satisfy the relation
\begin{equation} \label{eq:graph_assumption}
    \pi(\bfz \mid \bfx) = \pi(\bfz \mid \bfb).
\end{equation}
In this case, the ratio in \cref{eq:graph_assumption} is one, and $\pi(\bfx\mid\bfb,\bfz)$ coincides with the posterior predictive.

\subsection{Learning VED maps}
\label{sub:learnIPT}
In what follows, we discuss how to use training data to approximate maps $\widebar e$ and $\widebar d$ using VEDs. Assume that training data, $\{(\bfb^j, \bfx^j) \}_{j=1}^J$, containing observation and QoI pairs are provided. Let $e:\bbR^m \times \bbR^{p^{\rm e}}\to \bbR^{2\ell}$ and $d:\bbR^\ell\times \bbR^{p^{\rm d}} \to \bbR^q$ be an encoder and a decoder network, respectively.

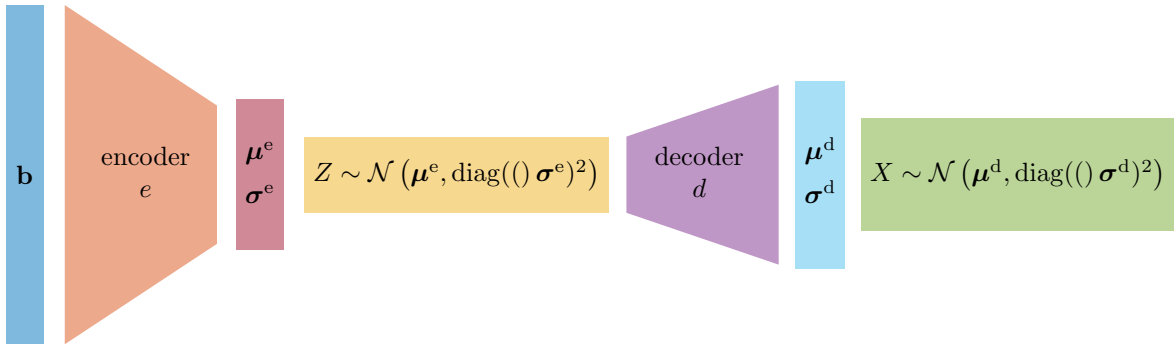
\begin{figure}[bthp]
    \centering
    \begin{tikzpicture}[scale = 0.92]
     
	\node[fill=matlab1!50, minimum width=0.5cm, minimum height=4.5cm] (X) at (0,0) {$\bfb$};
	
	\draw[fill=matlab2!50,draw=none] ([xshift=0.3cm]X.north east) -- ([xshift=2.5cm,yshift=1.0cm]X.east) -- ([xshift=2.5cm,yshift=-1.0cm]X.east) -- ([xshift=0.3cm]X.south east) -- cycle; 
	\node at (1.75,0) {$\begin{matrix}\mbox{encoder} \\ e \end{matrix}$};
	
	\node[fill=matlab7!50, minimum width=0.5cm, minimum height=2.0cm] (sa) at (3.4cm,0) {$\begin{matrix}\bfmu^{\rm e} \\[1ex] \bfsigma^{\rm e}\end{matrix}$};
	
	\node[fill=matlab3!50, minimum width=0.5cm, minimum height=1.0cm] (mu) at (6.24cm,0) {$Z \sim \calN\left(\bfmu^{\rm e},\diag(\bfsigma^{\rm e})^2\right)$};

	\draw[fill=matlab4!50,draw=none] ([xshift=0.25cm]mu.north east) -- ([xshift=2.45cm,yshift=0.75cm]mu.north east) -- ([xshift=2.45cm,yshift=-0.75cm]mu.south east) -- ([xshift=0.25cm]mu.south east) -- cycle;
	\node at (9.75,0) {$\begin{matrix}\mbox{decoder} \\ d \end{matrix}$};
	
	\node[fill=matlab6!50, minimum width=0.5cm, minimum height=2.5cm] (Xp) at (11.5,0) {$\begin{matrix}\bfmu^{\rm d} \\[1ex] \bfsigma^{\rm d}\end{matrix}$};
	
	\node[fill=matlab5!50, minimum width=0.5cm, minimum height=1.5cm] (XX) at (14.35cm,0) {$X \sim \calN\left(\bfmu^{\rm d},\diag(\bfsigma^{\rm d})^2\right)$};
     
\end{tikzpicture}
    \caption{Variational encoder-decoder network mapping for goal-oriented UQ. The encoder maps the input vector $\bfb$ to the latent variables $\bfz$ and the decoder maps the latent variables to the target $\bfx$.}
    \label{fig:VED_UQ}
\end{figure}

In \cref{fig:VED_UQ}, we provide a visual representation of the VED network used for goal-oriented UQ. The main distinction from the VED network in \Cref{fig:VEDnetwork} is that, in addition to enforcing a distribution on the latent variables, we assume a distribution on the target space so that we can perform UQ of the QoI.

Ideally, we would like to obtain the equality in \cref{eq:graph_assumption} with a learned encoder $e$ and decoder $d$. Hence, we aim to learn the parameterized encoder and decoder to minimize the distance between both distributions $\pi^{\rm e}(\bfz \mid \bfb)$ and $\pi^{\rm d}(\bfz \mid \bfx)$ in equation~\cref{eq:graph_assumption}.
We relax \cref{eq:graph_assumption} to obtain the minimization problem
\begin{equation} \label{eq:kl_minimization}
    \min_{\bftheta^{\rm e}, \bftheta^{\rm d}} \ \calD_{\rm kl}\left(\pi^{\rm e}(\bfz \mid \bfb) \mid\mid \pi^{\rm d}(\bfz \mid \bfx) \right).
\end{equation}
Here, $\pi^{\rm e}$ is a family of distributions that the encoder $e$ defines on the latent space. Similarly, $\pi^{\rm d}$ is a family of posterior distributions that are defined on the latent space. Therefore, the minimization problem \cref{eq:kl_minimization} looks for the best approximation to \cref{eq:graph_assumption} within these families. This approach is referred to as \emph{variational inference} \cite{goodfellow2016deep}. The superscript on a density function indicates the family to which the density function belongs.

We now expand the second argument in \cref{eq:kl_minimization} to obtain
\begin{equation}
    \begin{aligned}
    \calD_{\rm kl}\left(\pi^{\rm e}(\bfz \mid \bfb) \mid\mid \pi^{\rm d}(\bfz \mid \bfx) \right) &= \int_{\bbR^\ell} \pi^{\rm e}(\bfz \mid \bfb) \log \left( \frac{\pi^{\rm e}(\bfz \mid \bfb)}{\pi^{\rm d}(\bfz \mid \bfx)} \right) ~ \d\bfz \\
    &= \int_{\bbR^\ell} \pi^{\rm e}(\bfz \mid \bfb) \log \left( \frac{\pi^{\rm e}(\bfz \mid \bfb)\pi(\bfx)}{\pi^{\rm d}(\bfz,\bfx)} \right) ~ \d\bfz \\
    & =  \log \left( \pi(\bfx) \right)  \int_{\bbR^\ell} \pi^{\rm e}(\bfz \mid \bfb)~ \d\bfz + \int_{\bbR^\ell} \pi^{\rm e}(\bfz \mid \bfb) \log \left( \frac{\pi^{\rm e}(\bfz \mid \bfb)}{\pi^{\rm d}(\bfz,\bfx)} \right) ~ \d\bfz \\
    & = \log \left( \pi(\bfx) \right) - \calL(\bftheta^{\rm e}, \bftheta^{\rm d}).
    \end{aligned}
\end{equation}
Here we used the definition of a conditional density in the second line, and the logarithm product rule in the third line. Moreover, $\calL(\bftheta^{\rm e}, \bftheta^{\rm d})$ is called the \emph{evidence lower bound} (ELBO) \cite{goodfellow2016deep} and is defined as
\begin{equation} \label{eq:ELBO_definition}
    \calL(\bftheta^{\rm e}, \bftheta^{\rm d}) = \int_{\bbR^\ell} \pi^{\rm e}(\bfz \mid \bfb) \ \log \left( \frac{\pi^{\rm d}(\bfz,\bfx)}{\pi^{\rm e}(\bfz \mid \bfb)} \right) ~ \d\bfz.
\end{equation}
Note that $\log( \pi(\bfx) )$, also referred to as \emph{evidence}, does not depend on the maps $e$ and $d$. Therefore, it is a constant for minimization problem \cref{eq:kl_minimization}. We can drop this constant and reformulate \cref{eq:kl_minimization} to obtain
\begin{equation} \label{eq:elbo_minimization}
    \min_{\bftheta^{\rm e}, \bftheta^{\rm d}} \ - \calL(\bftheta^{\rm e}, \bftheta^{\rm d}).
\end{equation}
We now apply Bayes' theorem to \cref{eq:ELBO_definition} to obtain
\begin{equation} \label{eq:loss_final}
    \begin{aligned}
    \calL(\bftheta^{\rm e}, \bftheta^{\rm d}) &= \int_{\bbR^\ell} \pi^{\rm e}(\bfz \mid \bfb) \log \left( \frac{\pi^{\rm d}(\bfz,\bfx)}{\pi^{\rm e}(\bfz \mid \bfb)} \right) ~ \d\bfz \\
    &= \int_{\bbR^\ell} \pi^{\rm e}(\bfz \mid \bfb) \log \left( \frac{\pi^{\rm d}(\bfx \mid \bfz)\pi(\bfz)}{\pi^{\rm e}(\bfz \mid \bfb)} \right) ~ \d\bfz \\
    &= \int_{\bbR^\ell} \pi^{\rm e}(\bfz \mid \bfb) \log\left( \pi^{\rm d}(\bfx \mid \bfz)\right) ~ \d\bfz + \int_{\bbR^\ell} \pi^{\rm e}(\bfz \mid \bfb) \log \left( \frac{\pi(\bfz)}{\pi^{\rm e}(\bfz \mid \bfb)} \right) ~ \d\bfz \\
    & = \bbE_{\pi^{\rm e}(\bfz \mid \bfb)} \log \left( \pi^{\rm d}(\bfx \mid \bfz) \right) - \calD_{\rm kl}\left( \pi^{\rm e}(\bfz \mid \bfb)  \mid  \mid  \pi(\bfz)  \right).
    \end{aligned}
\end{equation}

The first term in the right-hand-side of \cref{eq:loss_final} is interpreted as computing an expectation with respect to the density function $\pi^{\rm e}(\bfz \mid \bfb)$. Note that the KL divergence in \cref{eq:kl_minimization} is with respect to the posterior density $\pi^{\rm d}(\bfz \mid \bfx)$ while the one in \cref{eq:loss_final} is with respect to the prior density $\pi(\bfz)$.

Computing the ELBO, as expressed in \cref{eq:loss_final}, still poses challenges, e.g., density functions $ \pi^{\rm e}(\bfz \mid \bfb)$ and $\pi^{\rm d}(\bfz  \mid \bfx)$ are intractable for nonlinear maps $e$ and $d$. We consider the following approximations to simplify the computation of the ELBO.

\begin{enumerate}
    \item[(a)] The prior density $\pi(\bfz)$ is a standard normal density, i.e., following the distribution $\mathcal N(\bfzero,\bfI)$, with $\bfI$ the identity matrix of size $\ell$.
    \item[(b)] We assume that $\pi^{\rm e}(\bfz \mid \bfb)$ is a product of independent Gaussian densities, i.e.,
    \begin{equation} \label{eq:latent-distribution}
        \pi^{\rm e}(\bfz  \mid \bfb) = \frac{1}{(2\pi)^{\ell/2} | \bfD^{\rm e}(\bfb)| } \exp \left(- \thf \| \bfD^{\rm e}(\bfb)^{-1} (\bfz - \bfmu^{\rm e}(\bfb)) \|_2^2 \right),
    \end{equation}
    where $\bfD^{\rm e}(\bfb) = {\rm diag} (\bfsigma^{\rm e}(\bfb))$ is a diagonal matrix and $|\cdot|$ denotes the determinant.
    \item[(c)] We assume that $\pi^{\rm d}(\bfx \mid \bfz)$ is a product of independent Gaussian densities of the form
 \begin{equation} \label{eq:output-distribution}
         \pi^{\rm d}(\bfx \mid \bfz) = \frac{1}{(2\pi)^{q/2} | \bfD^{\rm d}(\bfz)|} \exp\left( - \thf \| \bfD^{\rm d}(\bfz)^{-1} (\bfx - \bfmu^{\rm d}(\bfz)) \|_2^2 \right),
    \end{equation}
    where $\bfD^{\rm d}(\bfz) = {\rm diag} (\bfsigma^{\rm d}(\bfz)).$
\end{enumerate}

Combining approximations  (a) and (b), the KL divergence in \cref{eq:loss_final} can be written as
\begin{equation}
\calD_{\rm kl}\left( \pi^{\rm e}(\bfz \mid \bfb)  \mid  \mid  \pi(\bfz)  \right) =  \frac{1}{2}\left( -\log |\bfD^{\rm e}(\bfb)^2| - \ell + \bfmu^{\rm e}(\bfb) \t \bfmu^{\rm e}(\bfb) +
{\rm trace}\{ \bfD^{\rm e}(\bfb)^2 \}
\right).
\end{equation}
Then, with approximation (c), the ELBO simplifies to
\begin{equation} \label{eq:cost_final}
    \begin{aligned}
    \calL(\bftheta^{\rm e}, \bftheta^{\rm d}) = \bbE_{\pi^{\rm e}(\bfz \mid \bfb)}& \left(-\frac{q}{2}\log(2\pi) - \log | \bfD^{\rm d}(\bfz) |  - \frac{1}{2} \| \bfD^{\rm d}(\bfz)^{-1} (\bfx - \bfmu^{\rm d}(\bfz)) \|_2^2 \right) \\
    &- \frac{1}{2}\left( -\log |\bfD^{\rm e}(\bfb)^2| - \ell + \bfmu^{\rm e}(\bfb) \t \bfmu^{\rm e}(\bfb) +
{\rm trace}\{ \bfD^{\rm e}(\bfb)^2 \}
\right).
    \end{aligned}
\end{equation}

\begin{remark} In the case where $\bfD^{\rm d}(\bfz) \equiv \eta \bfI$, where $\eta > 0$ is a constant and $\bfI$ is an identity matrix, the cost function \cref{eq:cost_final} (sans constant) reduces to
\begin{equation} \label{eq:cost_final_simplified}
    \begin{aligned}
    \bbE_{\pi^{\rm e}(\bfz \mid \bfb)}& \left(- \frac{1}{2\eta^2} \| \bfx - \bfmu^{\rm d}(\bfz) \|_2^2\right) - \frac{1}{2}\left( -\log |\bfD^{\rm e}(\bfb)^2| - \ell + \bfmu^{\rm e}(\bfb) \t \bfmu^{\rm e}(\bfb) +
{\rm trace}\{ \bfD^{\rm e}(\bfb)^2 \}
\right).
    \end{aligned}
\end{equation}
A good guess for $\eta$ can significantly reduce the complexity of minimizing the cost function.
\end{remark}

Note that this is for one sample, so the overall goal is to minimize the total loss,
\begin{equation}
\label{eq:VED_problem}
\min_{\bftheta^{\rm e}, \bftheta^{\rm d}} \ \sum_{j=1}^J \calL^j(\bftheta^{\rm e}, \bftheta^{\rm d})
\end{equation}
where $\calL^j(\bftheta^{\rm e}, \bftheta^{\rm d})$ is defined \cref{eq:cost_final} for $\{\bfb^j, \bfx^j \}$. A standard optimization algorithm such as ADAM can be used to solve \cref{eq:VED_problem}.  See \cref{alg:VED_ADAM} for a general outline.

\begin{algorithm}[htb!]
\caption{Learning VED networks}
\label{alg:VED_ADAM}
    \begin{algorithmic}[1]
        \State \textbf{input} $\{(\bfb^j, \bfx^j)\}_{j=1}^J$, initialize encoder $e$ with $\bftheta^{\rm e}$, initialize decoder $d$ with $\bftheta^{\rm d}$
        \For{$k = 1, 2, \ldots$ (optimization steps) }
        \State Get minibatch
        \For{$j$ in the minibatch (loop over elements in minibatch)}
        \State Forward propagate encoder $[\bfmu^{\rm e},\bfsigma^{\rm e}] = e(\bfb^j;\,\bftheta^{\rm e})$
        \State Generate $L$ samples from the latent space
        $\bfz^l \sim \calN(\bfmu^{\rm e}, \diag{(\bfsigma^{\rm e})})$
             \For{$l=1, \ldots, L$  (loop over samples from latent space)}
            \State Evaluate decoder $[\bfmu^{\rm d},\bfsigma^{\rm d}] = d(\bfz^l;\,\bftheta^{\rm d})$
        \EndFor
        \State Compute loss $\calL^j$ and gradient $\nabla \calL^j$
        \EndFor
        \State Update $\bftheta^{\rm e}$ and $\bftheta^{\rm d}$
        \EndFor
        \State \textbf{output} $\bftheta^{\rm e}$ and $\bftheta^{\rm d}$
    \end{algorithmic}
\end{algorithm}

\subsection{Quantifying uncertainties in the QoI}
\label{sub:UQforQoI}
In this section, we discuss how to estimate and perform UQ on QoI $\bfx$ for a trained VED network, given an unobserved vector $\bfb$. In the context of goal-oriented inverse problems, a common approach to estimate the QoI is to use point estimates, e.g., the mean or the maximum of the posterior predictive. However, to quantify uncertainties for goal-oriented inverse problems one may need to compute the higher central moments of the posterior predictive, which in this case is computationally prohibitive. Hence, we first relate the outputs of the trained VED with the posterior predictive, describe how to sample from the posterior, and then utilize generated samples to estimate uncertainties for the QoI.

The following proposition relates the posterior-predictive with the output of the VED network.
\begin{proposition}
Let  $B,Z,$ and $X$ denote random variables with corresponding probability densities $\pi(\bfb),\pi(\bfz),$ and $\pi(\bfx)$. The posterior predictive is related to the output of the VED according to
\begin{equation} \label{eq:post-VED}
    \pi_{\rm pred}(\bfx \mid \bfb) = \mathbb E_{\bfz \mid \bfb}\ \pi(\bfx \mid \bfz, \bfb).
\end{equation}
\begin{proof}
\begin{equation} \label{eq:post-VED-proof1}
    \pi_{\rm pred}(\bfx \mid \bfb) = \int_{\mathbb R^{\ell}} \pi(\bfx,\bfz \mid \bfb)\ d\bfz = \int_{\mathbb R^{\ell}} \pi(\bfx \mid \bfb, \bfz) \pi(\bfz \mid \bfb) \ d\bfz = \mathbb E_{\bfz \mid \bfb}\ \pi(\bfx \mid \bfz, \bfb).
\end{equation}
Here, we used the definition of marginalization in the first equality and the conditional probability in the second equality. We refer to $\mathbb E_{\bfz \mid \bfb}\ \pi(\bfx \mid \bfz, \bfb)$ as the \emph{VED posterior predictive distribution}.

In the case where $\pi(\bfb\mid\bfz) = \pi(\bfb\mid\bfx,\bfz)$, as in \Cref{assump:bias}, we can further simplify \cref{eq:post-VED-proof1}. Let us first reformulate the joint distribution $\pi(\bfx,\bfz,\bfb)$.
\begin{equation}
    \pi(\bfx,\bfz,\bfb) = \pi(\bfb\mid\bfx,\bfz) \pi(\bfx,\bfz) = \pi(\bfb\mid\bfz) \pi(\bfx,\bfz) = \frac{\pi(\bfb,\bfz)}{\pi(\bfz)}\pi(\bfx,\bfz) = \pi(\bfb,\bfz)\pi(\bfx\mid\bfz).
\end{equation}
Here, we used the additional assumption in the second equality. Now we have
\begin{equation}
    \pi(\bfx \mid \bfb, \bfz) = \frac{\pi(\bfx,\bfz,\bfb)}{\pi(\bfz, \bfb)} = \frac{\pi(\bfb,\bfz)\pi(\bfx\mid\bfz)}{\pi(\bfb,\bfz )} = \pi(\bfx\mid\bfz).
\end{equation}
This simplifies \cref{eq:post-VED-proof1} to
\begin{equation} \label{eq:post-VED-proof2}
    \pi_{\rm pred}(\bfx \mid \bfb) = \mathbb E_{\bfz \mid \bfb}\ \pi(\bfx \mid \bfz).
\end{equation}

\end{proof}
\end{proposition}
Note that the expectation in the right-hand-side of \cref{eq:post-VED} is with respect to $\bfz$. Therefore, the VED-posterior is a probability distribution with respect to $\bfx$.

Now we discuss how to generate samples from the VED posterior-predictive and perform UQ, using the VED network. Suppose that an observation vector $\bfb$ is provided. We then form the random variable $Z\mid B$. Note that this has a Gaussian distribution $\mathcal N(\bfmu^{\text{e}},\diag{(\bfsigma^{\text{e}})}^2 )$. We generate samples $\{\bfz^l\}_{l=1, \ldots,L}$ from this distribution.  For each of these samples from the latent space, we evaluate the decoder $d(\bfz^l;\,\bftheta^{\rm d})$, and sample from the Gaussian distribution $\mathcal N(\bfmu^{\rm d},\diag{(\bfsigma^{\rm d})}^2 )$. The target  has $N_{\text{sample}}= L \cdot K$ samples from the VED posterior predictive.
\begin{algorithm}[htb!]
\caption{Sampling for go-UQ}\
    \begin{algorithmic}[1]
        \State \textbf{input} $\bfb$, encoder $e$ that depends on parameters $\bftheta^{\rm e}$, decoder $d$ that depends on parameters $\bftheta^{\rm d}$
        \State Evaluate encoder $[\bfmu^{\rm e},\bfsigma^{\rm e}] = e(\bfb;\,\bftheta^{\rm e})$
        \State Generate $L$ samples from the latent space
        $\bfz^l \sim \calN(\bfmu^{\rm e}, \diag{(\bfsigma^{\rm e})}^2)$
        \For{$l = 1, 2, \ldots, L$}
            \State Evaluate decoder $[\bfmu^{\rm d},\bfsigma^{\rm d}] = d(\bfz^l;\,\bftheta^{\rm d})$
            \State Generate $\kappa$ samples            $\bfx_k^l \sim \calN(\bfmu^{\rm d}, \diag{(\bfsigma^{\rm d})}^2)$
        \EndFor
        \State \textbf{output} $\{\bfx^l_k\}_{k=1, \ldots \kappa;\ l = 1, \ldots, L}$ samples from VED predictive posterior
    \end{algorithmic}
    \label{alg:samplingUQ}
\end{algorithm}

Note that we can estimate the mean of the VED-posterior $\bfmu_{\text{VED}}$ using
\begin{equation} \label{eq:VED-mean}
    \bfmu_{\text{VED}} = \mathbb E_{\bfx} \mathbb E_{\bfz\mid\bfb}  \pi(\bfx\mid\bfz, \bfb) \approx \frac{1}{L} \sum_{i \leq L} \mathbb E_{\bfz^l\mid\bfb} \pi(\bfx\mid\bfz^l, \bfb) = \frac{1}{L} \sum_{l \leq L} {\bfmu}^{\rm d}(\bfz^l).
\end{equation}

\section{Numerical Results}
\label{sec:numerics}
In this section, we consider two examples from inverse problems, where the aim is to perform UQ for some quantities of interest. The QoI is different in each example.  In \Cref{sub:experiment1}, the goal is to perform UQ on the regularization parameter required to compute tomography reconstructions using total variation regularization.  In \Cref{sub:experiment2}, we consider a hydraulic tomography example, where the QoI represents expansion coefficients defining a piecewise constant conductivity field.  In both instances, we sample from the predictive posterior and bypass the intermediate reconstructions (i.e., the tomography reconstruction and the conductivity field reconstruction, respectively).

\subsection{Experiment 1: X-ray Computed Tomography (CT)}
\label{sub:experiment1}
X-ray CT is an imaging technique that uses a combination of X-ray exposure and computational procedures to construct cross-sectional images of an object. When an object is exposed to an X-ray, its material attenuates the radiation whereas the amount of attenuation depends on the material's properties. The interaction between the object and the X-ray can be modeled as a line integral and is referred to as the Beer–Lambert law \cite{doi:10.1137/1.9781611976670}. The transformation that takes an object (an attenuation field) into its line integrals is mathematically modeled as a linear transformation and is referred to as the \emph{Radon transform} \cite{doi:10.1137/1.9781611976670}. The transformed function is referred to as a \emph{sinogram}. The CT reconstruction aims to infer the attenuation field from its Radon transform, however, this inversion process is ill-posed, and stabilization through regularization is required.

Let $\bfy \in \bbR^n$ be a vectorized 2-dimensional discretized and unknown attenuation field and let $\bfb \in \bbR^m$ be the noisy sinogram. Consider the model given in \cref{eq:detmodel}, where
$F(\bfy)$ denotes the discrete Radon transform of $\bfy$. The corresponding variational inverse problem can be written as a deterministic minimization problem
\begin{equation}\label{sub:tomoTV}
     \bfy(x) = \argmin_\bfy \ \norm[2]{F(\bfy) -\bfb}^2 + x \norm[1]{\bfD \bfy},
\end{equation}
where the first term represents the data fidelity and the second term is the regularization term.  Here we promote piecewise constant solutions in $\bfy$ by utilizing a 1-norm and a first-order two-dimensional finite difference operator $\bfD$, \cite{hastie2015statistical,rudin1992nonlinear}, representing a total variation regularizer. The balance between the two terms in this minimization problem is established by the regularization parameter $x >0$.

The regularization parameter is a priori unknown and choosing an appropriate $x$ can be a challenging task \cite{farquharson2004comparison,galatsanos1992methods,haber2000gcv,mead2008newton,vogel1996non}. Suppose that a true attenuation field $\bfy_{\text{true}}$ is available, then one way to select $x$ is to solve a bilevel optimization problem
\begin{equation} \label{eq:bilevel}
\widehat{x} =  \argmin_{x} \ \norm[2]{ \bfy(x) - \bfy_{\text{true}}}^2, \quad \mbox{where} \quad \bfy(x)  \quad \mbox{is given in \cref{sub:tomoTV}}.
\end{equation}
The optimization problem in \cref{sub:tomoTV} is regarded as the inner minimization problem with a fixed regularization parameter $x$. The minimization problem in \cref{eq:bilevel} is referred to as the \emph{design problem} which seeks the optimal regularization parameter. A bilevel optimization problem of type \cref{eq:bilevel} is, in general, non-convex.
The alternating direction method of multipliers (ADMM) method can be used to solve the inner-minimization problem and a Bayesian optimization method can be used to solve the outer-minimization problem \cite{chung2021efficient}. In \cite{afkham2021learning}, it was shown that DNNs can be trained to represent the mapping $\bfb \mapsto \widehat x$, but the framework did not allow for UQ.
Furthermore, it was observed that uncertainties in $\widehat x$ arise due to challenges in optimization in \cref{eq:bilevel} and numerical evaluation of its solution.
Thus, it is possible that a numerical solution for \cref{eq:bilevel} is far from the optimal regularization parameter.

\paragraph{Description of training} In this work, we train a VED network to represent the mapping $\bfb \mapsto \widehat x$.  Hence, in addition to obtaining a VED-predicted regularization parameter, we can efficiently draw samples from the VED predictive posterior and these samples allow us to quantify the uncertainty in the predicted parameter.  In this context, given $\bfy_{\rm true}$, the prediction operator $G$ corresponds to solving the bi-level minimization \cref{eq:bilevel} to find the optimal regularization parameter, i.e., $\widehat x = G(\bfy_{\rm true})$ is a numerically evaluated optimal regularization parameter. For the training data, we generate $J = 20,\!400$
images $\bfy_{\rm true}^{j}, j=1, \ldots, J$ as randomized Shepp-Logan phantoms \cite{randomSheppLogan,ruthotto2018optimal}. Then the training data are given by,
$
\left\{ F(\bfy_{\rm true}^{j}) + \bfe^{j}, G(\bfy_{\rm true}^{j})\right\}_{j=1}^J.$
For simplicity, we assume $F(\bfy) = \bfA \bfy$, where $\bfA$ is a matrix representing parallel beam tomography (with 181 parallel rays over 180 equidistant angles). The noisy sinogram is then given by
\begin{equation}
\label{eq:tomoforward}
\bfb^j = \bfA \bfy_{\rm true}^j + \bfe^j.
\end{equation}
The noise vector $\bfe^j$ in \cref{eq:tomoforward} is chosen to be white noise with a varying noise level that is selected uniformly random between 0.1\% and 5\%.
A noise level of $r$\% corresponds to a noise vector $\bfe^j$ satisfying $\norm[2]{\bfe^j}/\norm[2]{ \bfA \bfy_{\text{true}}^j } = r/100$.

\begin{figure}[bthp]
    \centering
    \includegraphics[width=0.8\textwidth]{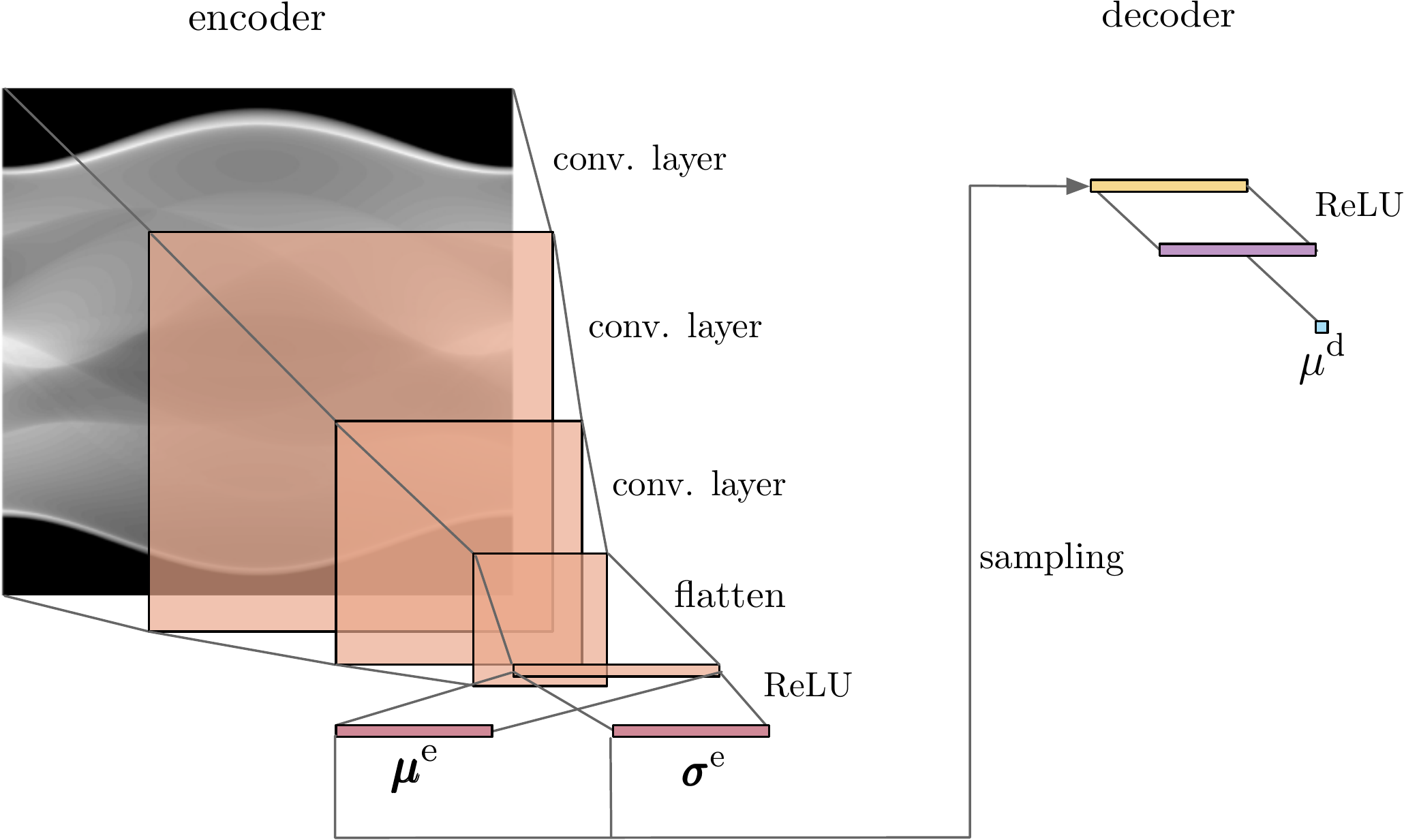}
    \caption{Architecture of the VED network for the X-ray CT problem, where the goal is to perform UQ on the regularization parameter.  On the left is the encoder network that maps the input sinogram to the latent space (defined by mean $\bfmu^{\rm e}$ and standard deviations $\bfsigma^{\rm e}$).  On the right is the decoder network that maps samples from the latent distribution to the output regularization parameter (defined by mean $\mu^{\rm d}$).}
    \label{fig:tomo_VED}
\end{figure}

The architecture of the VED network is illustrated in \cref{fig:tomo_VED}. The encoder network is comprised of convolutional layers and ReLU layers \cite{goodfellow2020generative}. The dimension of the latent space is set to $\ell=128$. The decoder network is a feed-forward fully-connected ReLU network with one hidden layer. This architecture is chosen to exploit local correlations in the intensities in a sinogram, but other architectures can be considered. We observe that other network architectures result in similar results (data not shown). The VED is trained by minimizing the cost function \cref{eq:cost_final_simplified}, with $\eta = 0.1$ using the ADAM optimization method \cite{kingma2014adam}. The optimization process is stopped after $10^4$ epochs. The learning rate is adjusted dynamically from the interval $[5\cdot 10^{-2},10^{-4}]$.

\begin{figure}
\begin{center}
     \begin{tabular}{ccc}
        \includegraphics[width=0.22\textwidth]{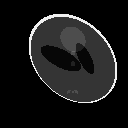} &
        \includegraphics[width=0.22\textwidth]{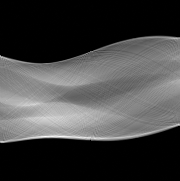} &
        \includegraphics[width=0.45\textwidth]{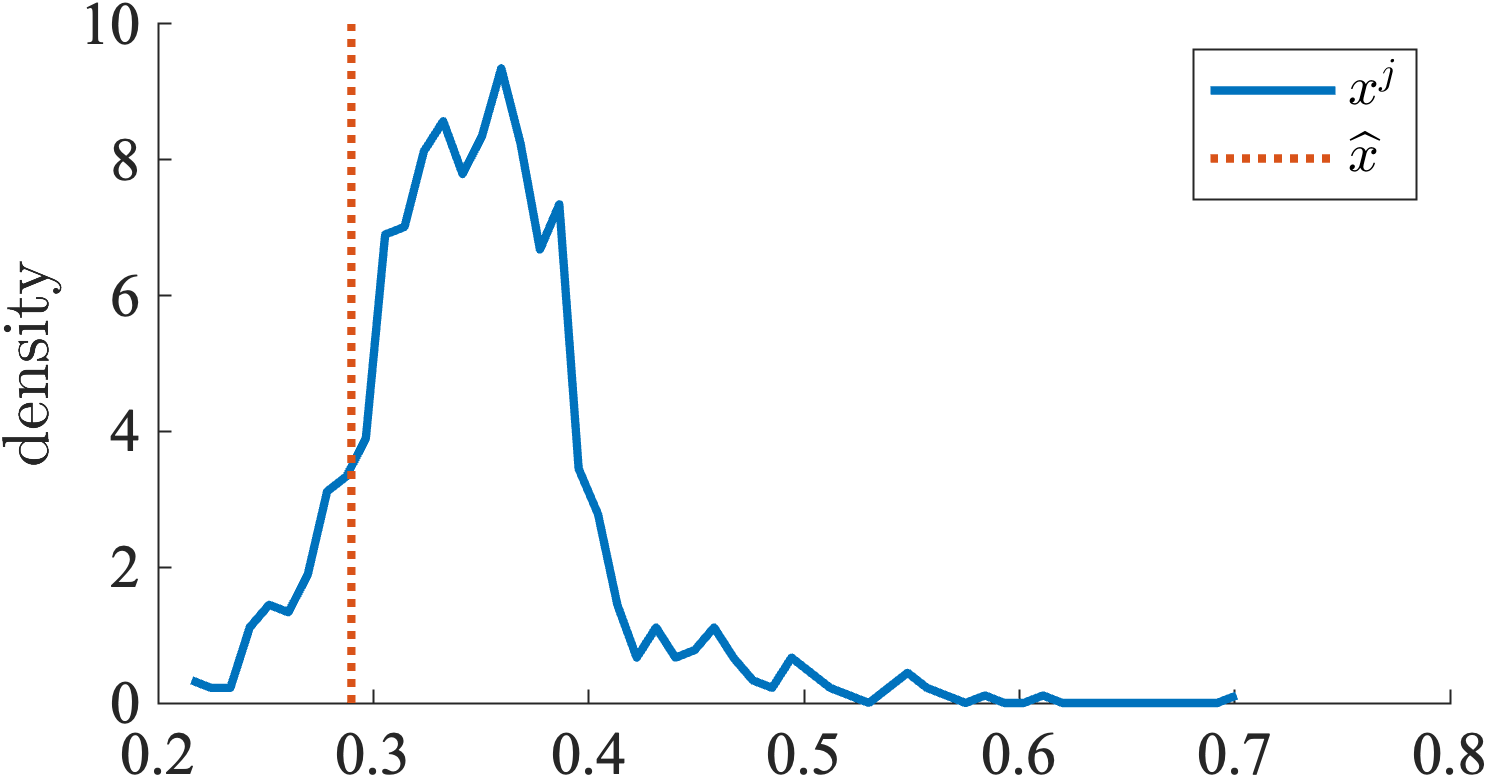}
        \\
        true image $\bfy_{\rm true}$  & sinogram $\bfb$ & distribution $\pi_{\rm pred}(x \mid \bfb)$ \\      \\
    \end{tabular}
     \begin{tabular}{ccc}
     \includegraphics[width=0.22\textwidth]{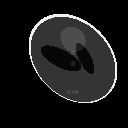} &
        \includegraphics[width=0.22\textwidth]{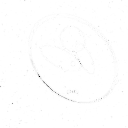} &
        \includegraphics[width=0.45\textwidth]{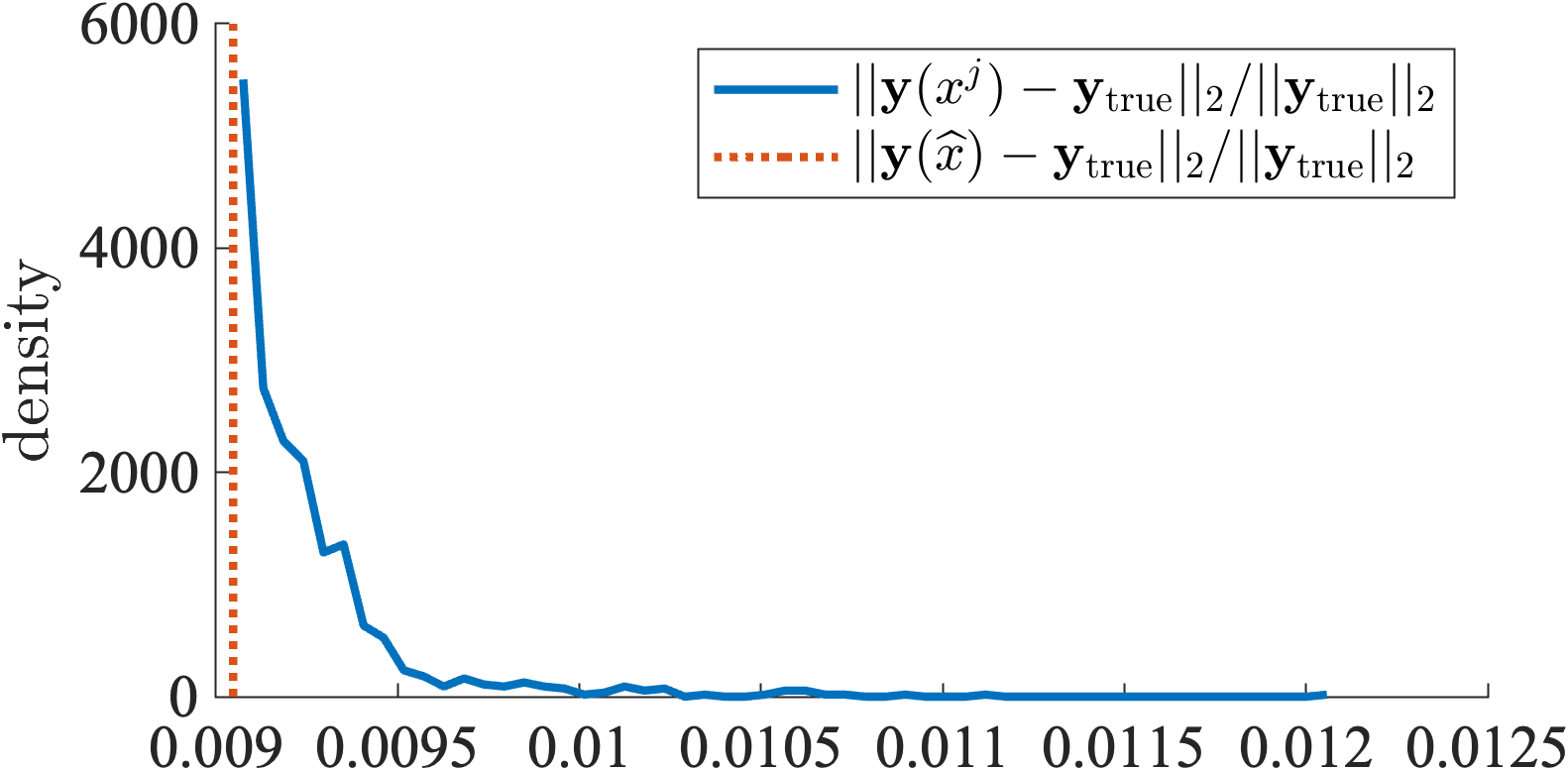}\\
        mean reconstruction &  pixel-wise variance &
    distribution rel. reconstruction errors \\
    \end{tabular}
\end{center}
    \caption{The trained VED network from the X-ray CT experiment is used to perform UQ on the regularization parameter for a test data sample.
     In the top row, we provide the true image, the corresponding noisy sinogram, and the distribution of the predicted regularization parameters $x^j$ from $\pi_{\rm pred}(x \mid \bfb)$ with the QoI $\widehat x$ provided for reference. In the bottom row, we provide the  mean reconstruction, the image of pixel-wise uncertainties (in inverted colormap so that large variances are denoted in black), and the distribution of the relative reconstruction errors. }
    \label{fig:tomography_test_data}
\end{figure}

\paragraph{Evaluation of VED}
To test the performance and generalizability of the learned VED network, we consider three datasets outside the training dataset. First, we consider an example from the test data. We generate a sinogram from a randomized Shepp-Logan phantom and add noise. The true image and sinogram are provided in the top row of \cref{fig:tomography_test_data}.  Using the approach described in \Cref{sub:UQforQoI}, we generate $10^3$ samples $x^j$ from the VED predictive posterior.  The distribution of these samples is provided in the top row of \cref{fig:tomography_test_data}.   Notice that the VED predicted values for the regularization parameter are within the interval $[0.2, .7]$, with a high probability of being close to the optimal regularization parameter $\widehat x = G(\bfy_{\rm true})$.
 To visualize the corresponding uncertainties in the reconstructed images, for each of the samples of the regularization parameter, we solve \cref{sub:tomoTV} to get a reconstruction $\bfy(x^j)$.
The mean reconstruction image, an image of the pixel-wise variances (in inverted colormap so that white corresponds to small variances and black corresponds to larger variances), and the distribution of relative reconstruction errors are provided in the bottom row of \cref{fig:tomography_test_data}.  Similar to observations made in \cite{uribe2022hybrid,afkham2023uncertainty}, we observe that large uncertainties occur at the locations of discontinuities in the image.  Furthermore, as expected, the relative reconstruction errors for samples from the posterior predictive are larger than the relative reconstruction error corresponding to the optimal regularization parameter (denoted by the red dotted vertical line), but the distribution of relative reconstruction errors is tight, implying a high probability of a good reconstruction.

\begin{figure}
\begin{center}
     \begin{tabular}{ccc}
        \includegraphics[width=0.22\textwidth]{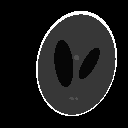} &
        \includegraphics[width=0.22\textwidth]{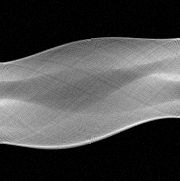} & \includegraphics[width=0.45\textwidth]{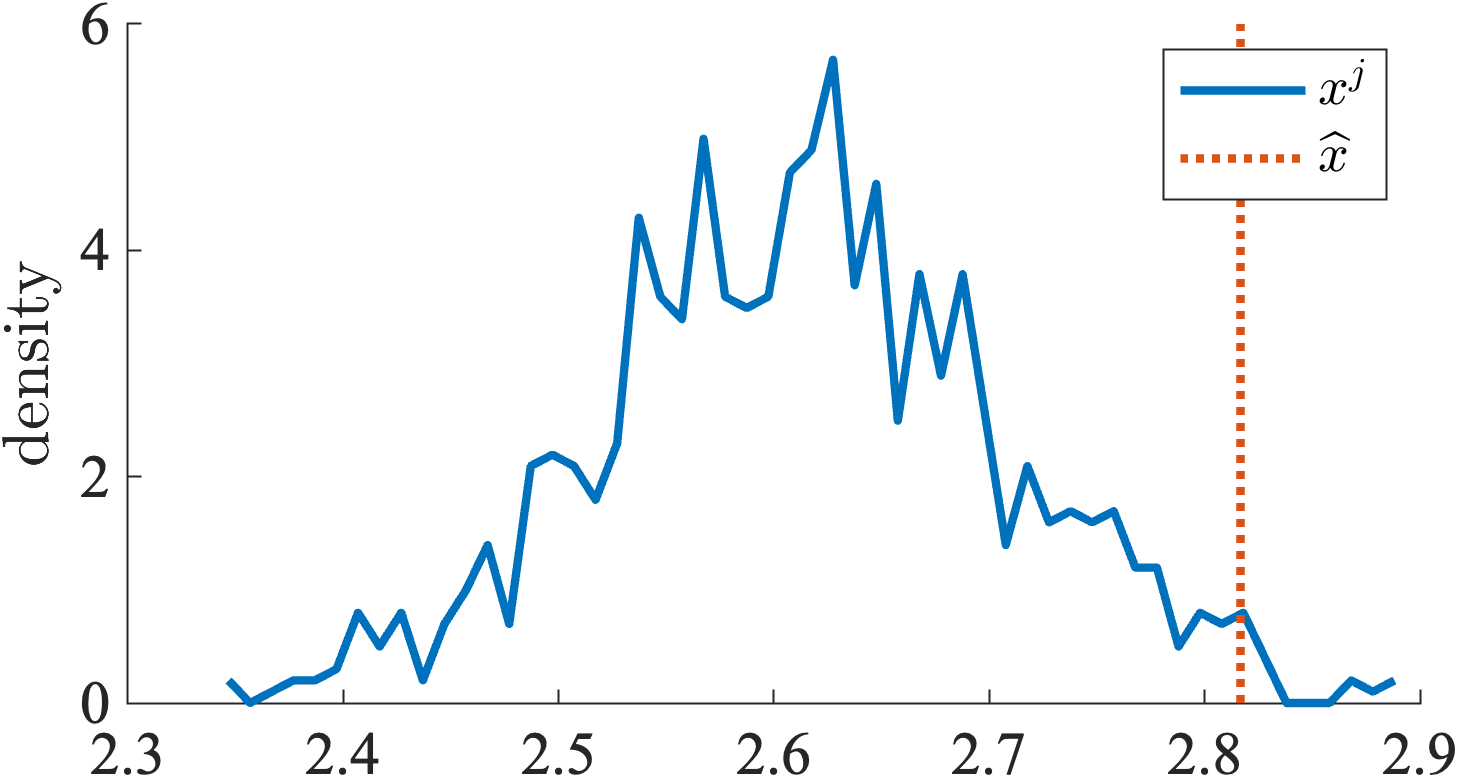}\\
        \\
        true image $\bfy_{\rm true}$  & sinogram $\bfb$ & distribution $\pi_{\rm pred}(x \mid \bfb)$  \\
    \end{tabular}
     \begin{tabular}{ccc}

        \includegraphics[width=0.22\textwidth]{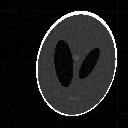} &
        \includegraphics[width=0.22\textwidth]{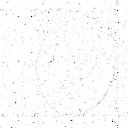}
       &
        \includegraphics[width=0.45\textwidth]{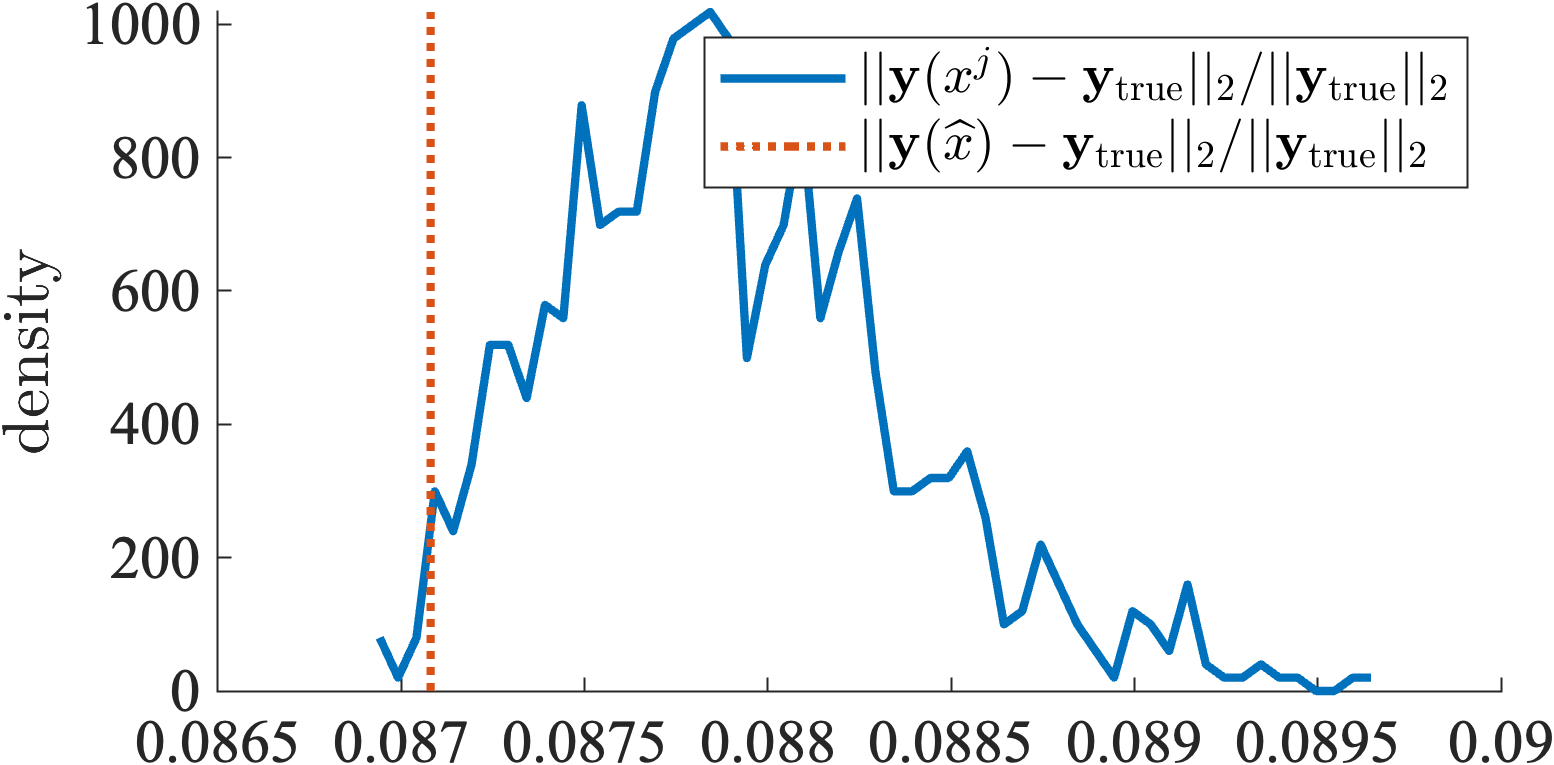}\\
       mean reconstruction &  pixel-wise variance  & distribution rel. reconstruction errors\\
    \end{tabular}
\end{center}
    \caption{UQ on the regularization parameter via a trained VED network for x-ray CT reconstruction with modified angles.}
    \label{fig:tomography_slightly_modified_A}
\end{figure}

For the second out-of-sample dataset, we generate a sinogram from a randomized Shepp-Logan phantom, but we modify the projection angles. That is, we use a forward model matrix $\bfA$ where the $180$ projections angles are no longer equidistant.  We add white noise to the vector of projection angles with a standard deviation of $0.1$.
The sinogram is propagated through the trained VED network, and we draw $10^3$ samples of the regularization parameter from the VED-posterior predictive $\pi_{\rm pred}(x \mid \bfb)$. The true image, the sinogram corresponding to modified angles, and the distribution of the regularization parameter samples are provided in the top row of \cref{fig:tomography_slightly_modified_A}.
Although the QoI $\widehat x$ is slightly larger than many of the predicted values, the range is fairly small.
The mean reconstruction, pixel-wise variances for the reconstruction, and the distribution of relative reconstruction errors are provided in the bottom row of \cref{fig:tomography_slightly_modified_A}. Obtaining relative errors smaller than the optimal can be attributed to numerical errors in solving \cref{eq:bilevel}.

\begin{figure}
\begin{center}
     \begin{tabular}{ccc}
        \includegraphics[width=0.22\textwidth]{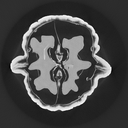} &
        \includegraphics[width=0.22\textwidth]{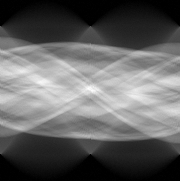} &
        \includegraphics[width=0.45\textwidth]{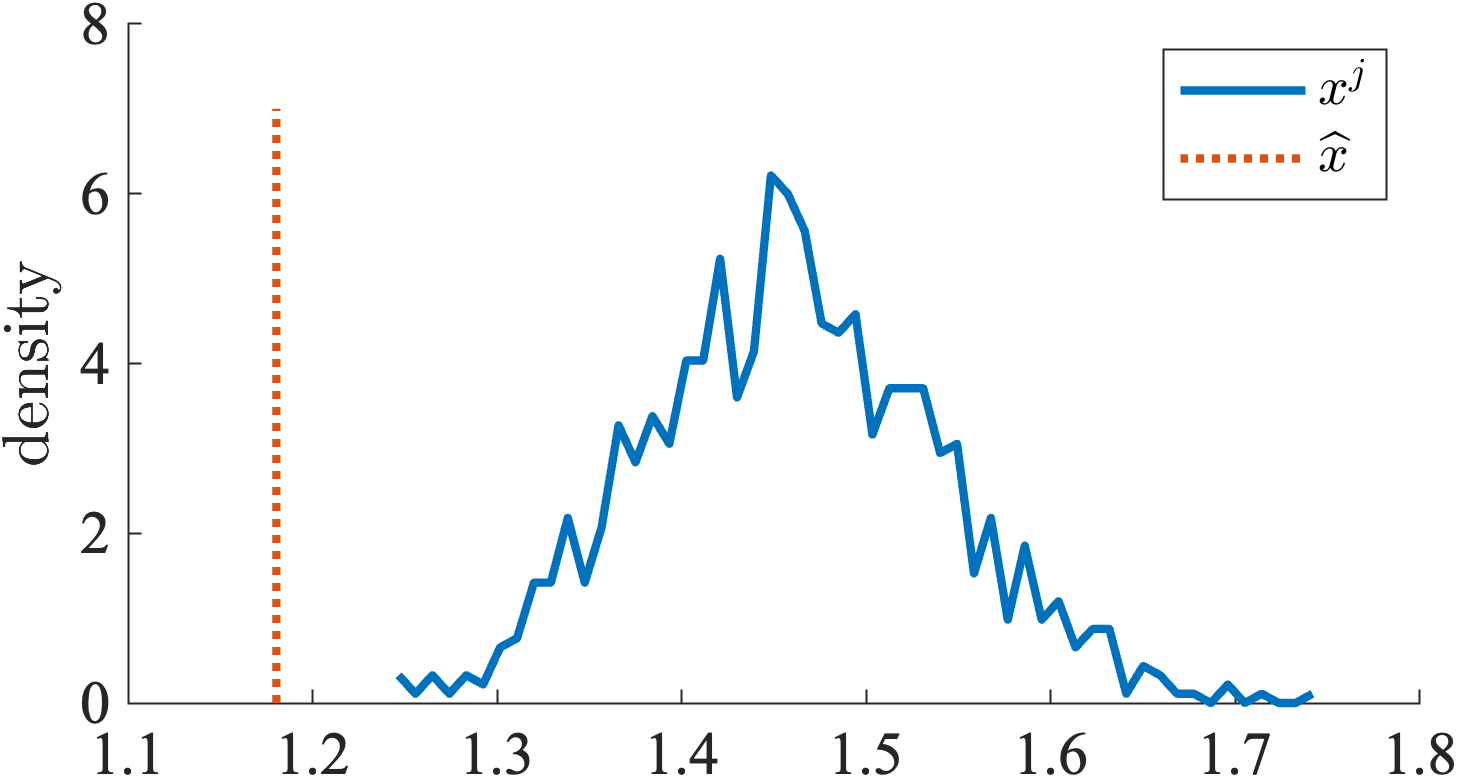}
        \\
        true image $\bfy_{\rm true}$  & sinogram $\bfb$ &  distribution $\pi_{\rm pred}(x \mid \bfb)$\\
    \end{tabular}
     \begin{tabular}{ccc}

        \includegraphics[width=0.22\textwidth]{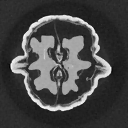} &
        \includegraphics[width=0.22\textwidth]{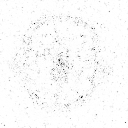} &
        \includegraphics[width=0.45\textwidth]{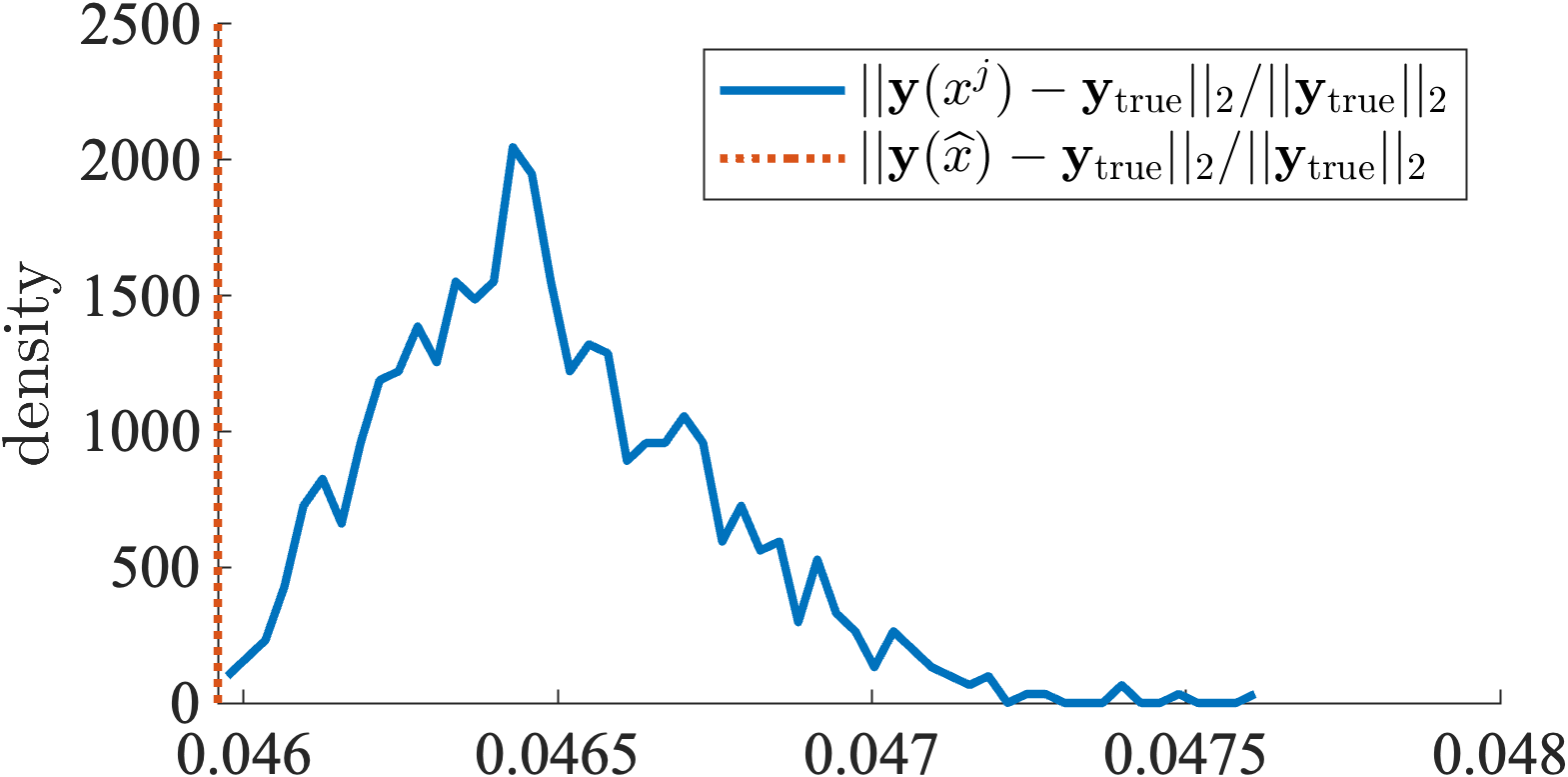}\\
        mean reconstruction &  pixel-wise variance  &  distribution rel. reconstruction errors \\
    \end{tabular}
\end{center}
    \caption{UQ on the regularization parameter via a trained VED network for an out-of-sample walnut image, \cite{https://doi.org/10.5281/zenodo.1254206}.}
    \label{fig:tomography_different_image}
\end{figure}

For the third out-of-sample experiment, we use data for the tomographic reconstruction of a walnut given at \cite{https://doi.org/10.5281/zenodo.1254206}. The walnut image, which was computed using a CT scan but that we take as ground truth, is provided in the top row of \cref{fig:tomography_different_image} along with a noisy sinogram that was generated as in \cref{eq:tomoforward}.
Similar to the previous experiment, we draw $10^3$ samples of the regularization parameter from the VED-posterior predictive $\pi_{\rm pred}(x \mid \bfb)$.  For this example, we observe that the trained VED network provides a wide range of values for the predicted regularization parameter, but all of the samples are not too far from the QoI $\widehat x$.  The mean reconstruction and pixel-wise variances are provided in the bottom row of \cref{fig:tomography_different_image}, along with the distribution of relative reconstruction errors.
Again, we see that large uncertainties are present at the locations of the discontinuities in the image, and the relative reconstruction error norms for the sampled regularization parameters are close to that of the optimal regularization parameter for this problem.

\subsection{Experiment 2: Nonlinear Hydraulic Tomography}
\label{sub:experiment2}
Next, we consider an example where the goal is to detect the locations of discontinuities in reconstructed  images.  A typical approach (e.g., used in CT) is to obtain an image reconstruction using standard algorithms, followed by post-processing, e.g., with image segmentation.  We consider an approach that avoids the intermediate step and goes directly from the observation to the discontinuities \cite{dahl2017computing,afkham2023uncertainty}.  Although the QoI could represent locations of the discontinuities, instead we consider a nonlinear hydraulic tomography problem that is adopted from \cite{cardiff2009bayesian,lee2013bayesian,2111.15620} and define the QoI to be coefficients representing a parameterization of the locations of the discontinuities.  Then, we use VEDs to obtain uncertainties for the coefficients.

Consider a confined aquifer that is modeled using an elliptic partial differential equation (PDE),
\begin{equation} \label{eq:aquifer}
    \begin{aligned}
    \nabla \cdot ( y(\bfxi) \nabla u(\bfxi) ) = q_i\delta(\bfxi^{\text{well}}_i),  \qquad \bfxi = (\xi_1,\xi_2)\in \calD, \\
    \end{aligned}
\end{equation}
where $\calD \subset \bbR^2$ is the unit square domain, $u:\mathbb R^2\to \mathbb R$ is the hydraulic head, $y:\mathbb R^2\to \mathbb R^+$ is the hydraulic conductivity, $\{\bfxi^{\text{well}}_i\}_{i=1}^{N_{\text{well}}}$ are locations for $N_{\text{well}}$ wells, and, $q_i$ is the pumping rate at the $i$-th well. We consider a zero-head (homogeneous Dirichlet) boundary condition on the left, right, and bottom boundaries of $\calD$. Furthermore, we consider a no-flux (homogeneous Neumann) boundary condition on the top boundary. We discretize and solve \cref{eq:aquifer} using a finite element method (FEM) on a structured mesh with $N_{\text{FEM}} = 10,201$ degrees of freedom (on a grid of approximate size $100\times100$).

\begin{figure}
     \centering
     \begin{subfigure}[b]{0.45\textwidth}
         \centering
    \includegraphics[width=0.7\textwidth]{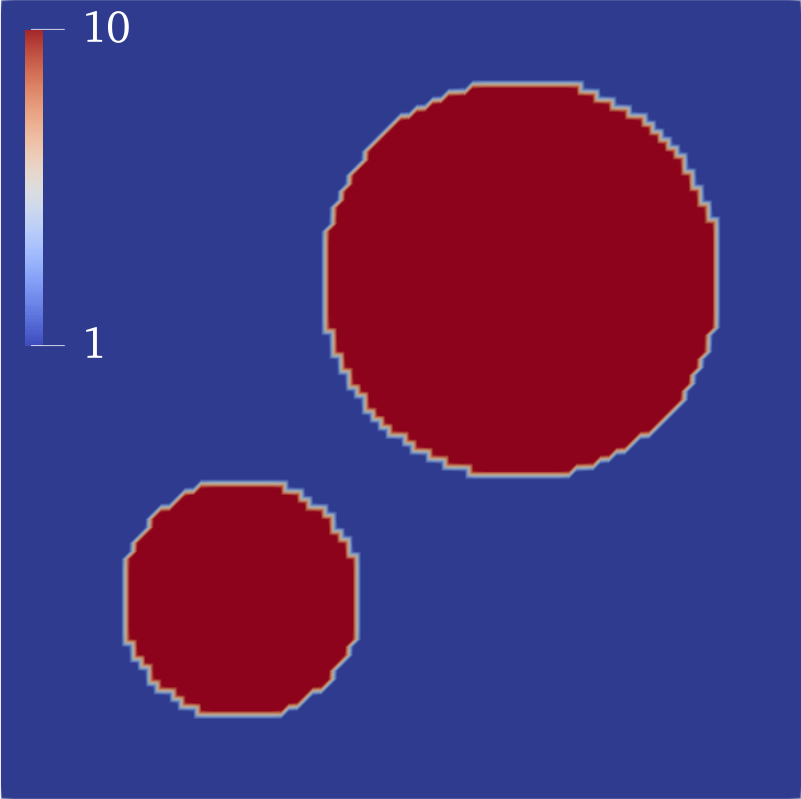}
         \caption{true hydraulic conductivity}
         \label{fig:true_conductivity}
     \end{subfigure}
     \begin{subfigure}[b]{0.45\textwidth}
         \centering
         \includegraphics[width=\textwidth]{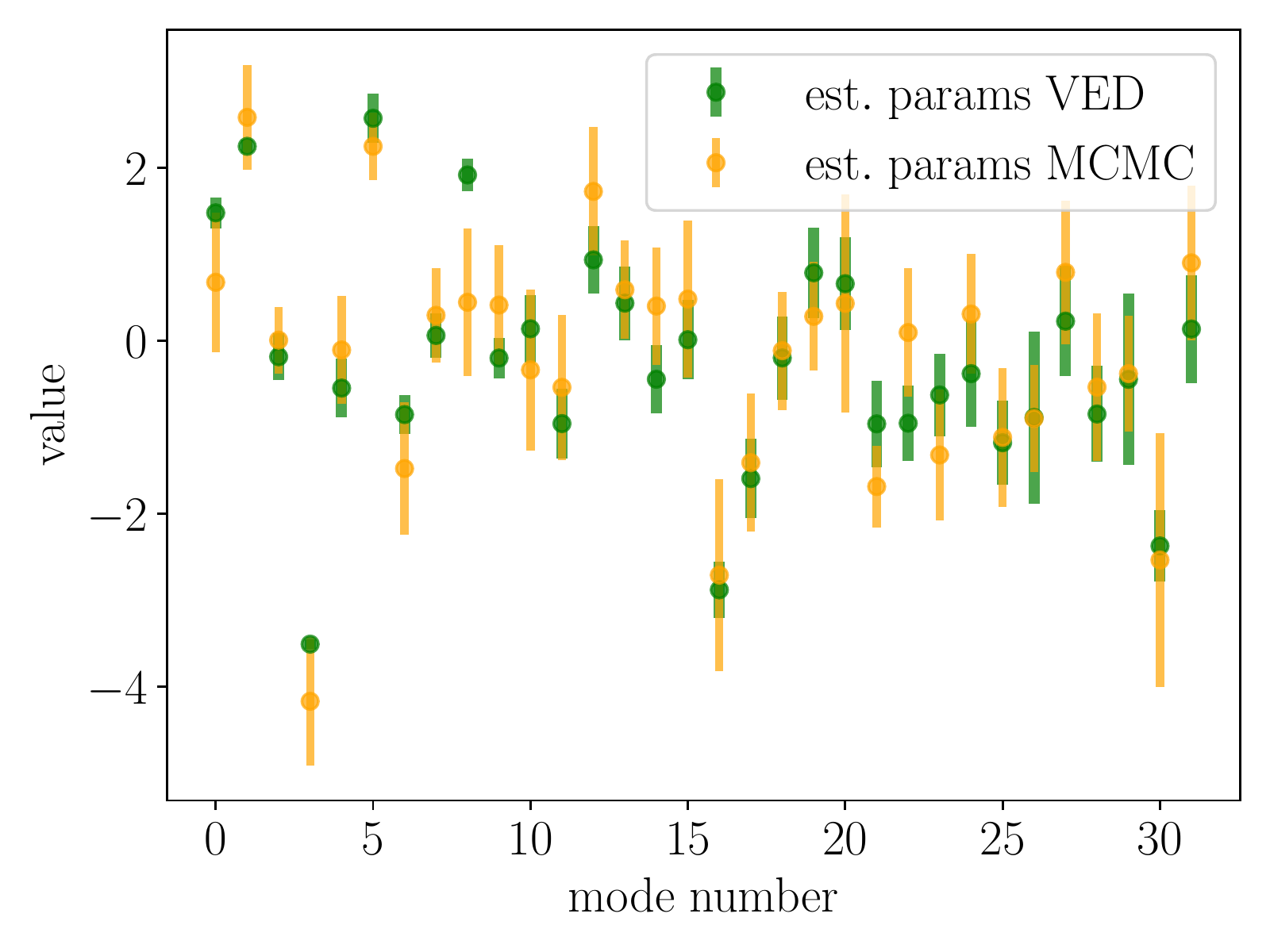}
         \caption{estimation of the expansion coefficients }
         \label{fig:est_expansion}
     \end{subfigure}
    \caption{Results for the nonlinear hydraulic tomography example. In (a), we provide the true hydraulic conductivity field.  In (b), we provide the estimated expansion coefficients and associated 99\% credibility intervals (indicated by bars) using samples from the VED posterior predictive and using MCMC samples. The mode number corresponds to the index $i$ for coefficient $X_i$.}
    \label{fig:hydraulic_diagnostics}
\end{figure}

We begin by describing the forward process. Assume we have the domain $\calD$ and the true hydraulic conductivity provided in \cref{fig:true_conductivity}. Consider placing 20 wells in the domain $\mathcal D$ located at coordinates $(i/4, j/5)$, with integers $1 \leq i \leq 4$ and $1 \leq j \leq 5$.  At the $i$-th well, we perform an injection at the rate of $q_i$. We measure the head at the rest of the $(N_{\text{well}}-1)$ wells and repeat the injection for all wells.  Thus, our measurement vector $\bfb$ is of size $N_{\text{well}}(N_{\text{well}} -1)$, where the $k$-th component of $\bfb$ is given as
\begin{equation} \label{eq:hydraulic_forward}
    b_k = u^{i}_{y}(\bfxi_j) + e_k, \quad k = i\cdot N_{\text{well}} + j, \text{ and, } i\neq j,
\end{equation}
where $u_y^{i}$ is the solution to \cref{eq:aquifer} when the conductivity is $y$ and the injection is at the $i$-th well, and $e_k$ represents error in the observation.
Note that the parameter to observation map $F$ in \cref{eq:detmodel} corresponds to the map $y\mapsto \bfb$. We consider a noise level of $1\%$, i.e., $\bfe$ containing errors $e_k$ is a realization of a mean-zero Gaussian multivariate random variable with covariance matrix $\sigma_n^2 \bfI$ where
\begin{equation}
    \sigma_n = \frac{\sqrt{ \sum_{i=1}^{N_{\text{well}}} (u_y^{i})^2  }}{100}.
\end{equation}
Given measurements $\bfb$, the aim of the hydraulic tomography inverse problem is to reconstruct the hydraulic conductivity.

We define a goal-oriented prior distribution for $y(\bfxi)$. Let $V$ be a zero-mean multivariate Gaussian random vector of size $N_{\text{FEM}}$, distributed according to the covariance operator $\mathcal C = (\epsilon^2_{\text{cl}} \bfI + \Delta)^{-1}$, with $\epsilon_{\text{cl}} = 1/0.1$, and, $\bfI$ and $\bfL$, the identity and the Laplacian operators, respectively. We now construct a piece-wise constant conductivity field
\begin{equation} \label{eq:prior_conductivity}
    y(\bfxi) = a^+ H( V(\bfxi) ) + a^-( 1-H( V(\bfxi) ) ),
\end{equation}
where $H$ is the Heaviside step function, and, $a^+, a^- > 0$ are positive constants. Note that we can approximate $V$ using the eigendecomposition of $\mathcal C$ as
\begin{equation} \label{eq:KL}
    V(\bfxi) \approx \sum_i^{N_{\text{goal}}} X_i \sqrt{g_i} v_i(\bfxi),
\end{equation}
where $\{ (g_i, v_i(\bfxi)) \}_{i \leq N_{\text{goal}}}$ are the eigenpairs of $\mathcal C$ associated with the largest eigenvalues of $\mathcal C$ sorted in decreasing order. Furthermore, $X_i$'s are independent and standard normal random variables. We define the goal of the inverse problem to be the coefficients $X_i$ in \cref{eq:KL}.
Let $N_{\text{goal}}=32$, then the QoI are parameters $\bfx \in \bbR^{32}$, and the prediction operator $G$, in \cref{eq:detmodel} corresponds to the map $y\mapsto \bfx$.

\begin{figure}[bthp]
    \centering
    \includegraphics[width=0.8\textwidth]{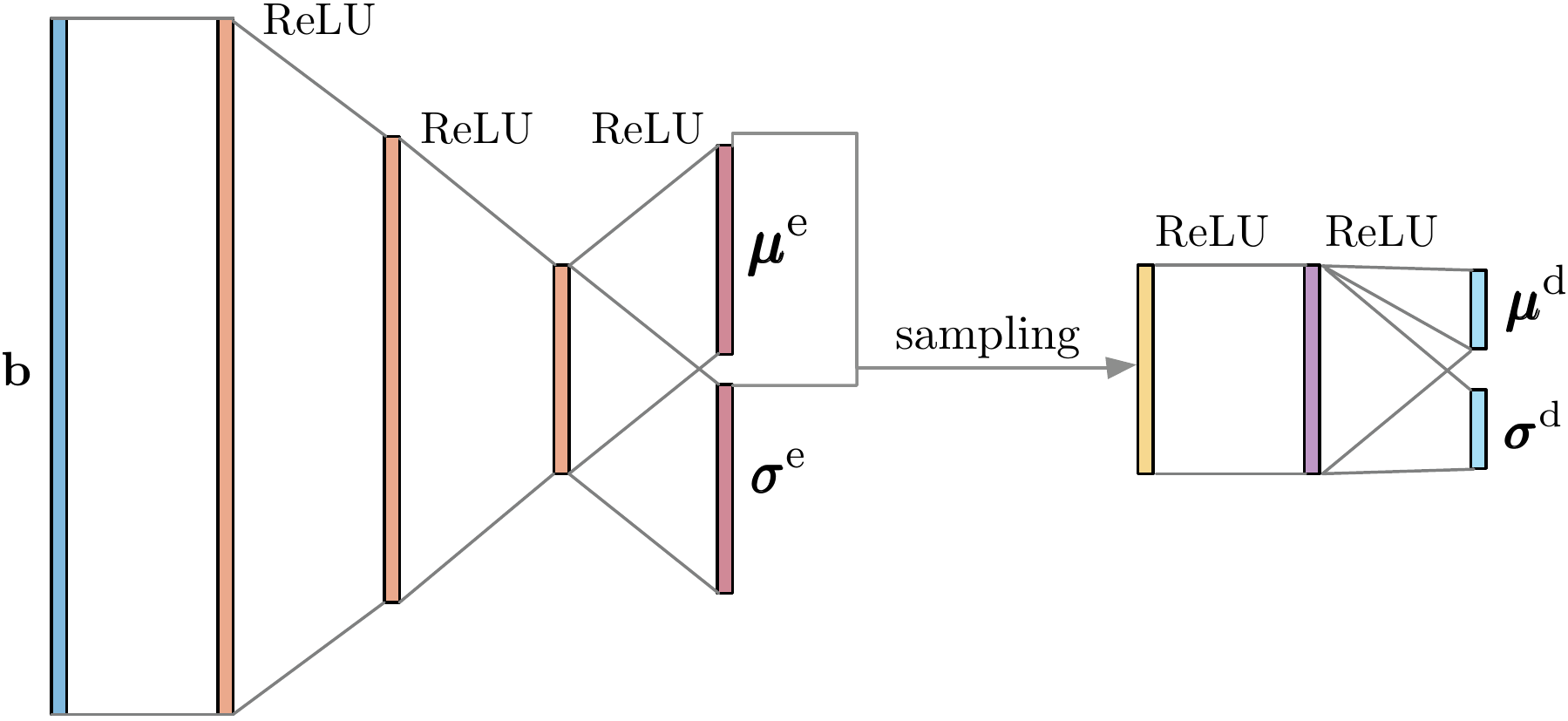}
    \caption{Architecture of the VED network for the hydraulic tomography problem, where the goal is to perform UQ on the parameters describing the discontinuities in the conductivity. On the left is the encoder network that maps the input, hydraulic head at the location of the wells, to the latent space (defined by mean $\bfmu^{\rm e}$ and standard deviations $\bfsigma^{\rm e}$).  On the right is the decoder network that maps samples from the latent distribution to the target discontinuity parameters (defined by mean $\bfmu^{\rm d}$ and standard deviations $\bfsigma^{\rm d}$).}
    \label{fig:hydraulic_VED}
\end{figure}

\paragraph{Description of training} We collect a dataset of pairs $\{ (\bfb^j,\bfx^j) \}_{j\leq J}$, with $J = 10^4$ where $\bfx^j$ are realizations from the standard normal distribution, and $\bfb^j$ has components described in \cref{eq:hydraulic_forward}.  We then consider the VED network summarized in \cref{fig:hydraulic_VED} for the hydraulic tomography problem. Note that in this example, we assume that we do not know the mean and the standard deviations of the target. We consider a Gaussian multi-layer perceptron (MLP) output layer that estimates these values with 2 output layers.
The cost function for this network can be constructed following \cref{eq:cost_final}, and the network is trained using the ADAM optimization algorithm \cite{kingma2014adam} over $1.5\times 10^4$ epochs. The learning rate is dynamically adjusted from the interval $[5\times 10^{-2}, 10^{-4}]$.

\paragraph{Evaluation of VED}
We evaluate the performance of the trained VED network for an out-of-sample test problem of finding the hydraulic conductivity field in \cref{fig:true_conductivity}, which is constructed by placing disks centered at coordinates $(0.3,0.65)$ and $(0.25,0.65)$ with radii $0.15$ and $0.25$ and conductivity values set to $a^+=10$ and $a^-=1$, inside and outside of the disk inclusions, respectively. Note that this conductivity field does not belong to the training dataset and is not generated from the prior distribution \cref{eq:KL} and \cref{eq:prior_conductivity}. A measurement vector $\bfb$ is generated as in \cref{eq:hydraulic_forward}, and we follow \cref{alg:samplingUQ} to draw samples from the VED posterior predictive distribution.  In \cref{fig:est_expansion}, we provide the mean of the posterior-predictive, i.e., the $32$ expansion coefficients in $\bfx$ using the VED (in green) along with the 99\% credibility intervals. For comparison purposes, we use a Markov chain Monte Carlo (MCMC) method to obtain a sequence of samples $\{\mathbf  x_{\text{MCMC}}^{j} \}_{j=1}^{N_{\text{MCMC}}}$, from the posterior-predictive. MCMC methods allow for estimation of moments, e.g., $\mathbb E(f(X))$ can be estimated using the ergotic sum
\begin{equation} \label{eq:ergotic}
    \mathbb E(f(X)) \approx \frac{1}{N_{\text{MCMC}}} \sum_{j=1}^{N_{\text{MCMC}}} f(\mathbf x_{\text{MCMC}}^{j}).
\end{equation}
We estimate the mean of the posterior-predictive $\mathbb E(X)$ by taking $f$ to be the identity map in \cref{eq:ergotic}. Similarly, we estimate the variance by taking $f(X) = \left(X - \mathbb E(X) \right)^2 $.
The preconditioned Crank-Nicolson (pCN) method \cite{10.1214/13-STS421} was used for this problem because of its suitability for expansions of the form \cref{eq:KL} \cite{10.1214/13-STS421}.
We provide means and variances of the estimated expansion coefficients (in yellow) in \cref{fig:est_expansion}.
We notice a correlation between the two methods, suggesting that the mean of the VED-posterior predictive correctly approximates the mean of the actual posterior predictive.

\begin{figure}
     \centering
     \begin{subfigure}[b]{0.45\textwidth}
         \centering
         \includegraphics[width=\textwidth]{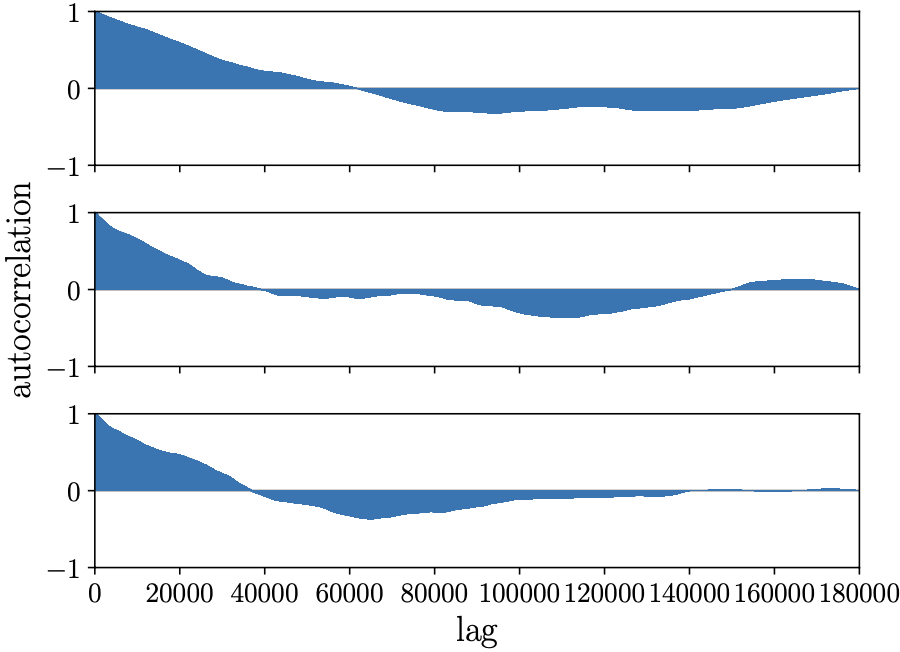}
         \caption{auto-correlation function of the first 3 coefficients.}
         \label{fig:acf}
     \end{subfigure}
     \begin{subfigure}[b]{0.45\textwidth}
         \centering
         \includegraphics[width=\textwidth]{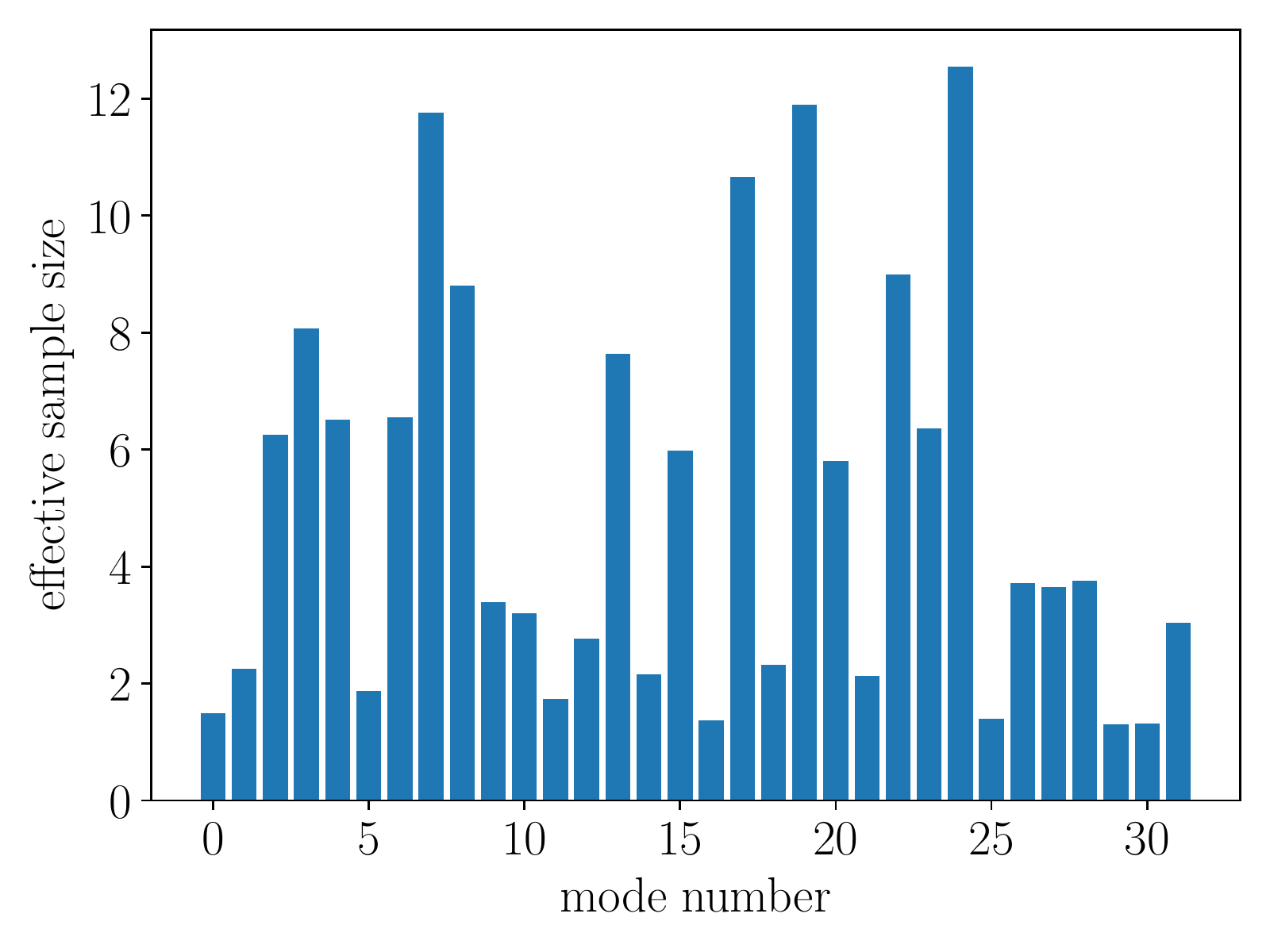}
         \caption{effective sample size for each expansion coefficient}
         \label{fig:ess}
     \end{subfigure}
    \caption{For MCMC samples from the posterior predictive, we provide the autocorrelation function for the first three expansion coefficients in (a) and the effective sample size for each coefficient in (b).}
    \label{fig:hydraulic_diagnostics2}
\end{figure}

The convergence of the sum in \cref{eq:ergotic} for a general posterior distribution is an ongoing topic of research. One indication for convergence is the effective sample size (ESS) \cite{mcbook} of MCMC samples, which suggests how many of the MCMC samples can be regarded as independent samples. Larger ESS values indicate a better approximation in \cref{eq:ergotic}.
In this example, we compute $2\times 10^5$ samples from the posterior predictive using the pCN sampler and discard the first $50\times 10^4$ samples as the burn-in period. \cref{fig:acf} shows the autocorrelation function \cite{mcbook} of the samples for the first 3 expansion coefficients in \cref{eq:KL}. We observe a similar behavior for the other coefficients. We notice that it takes about $50\times 10^4$ samples for the autocorrelation function to reach zero. This suggests that the convergence of the MCMC is only achieved with a very large number of samples. In \cref{fig:ess}, we provide the ESS for each expansion coefficient.
A reliable estimate of the variance of the posterior requires an ESS of around 100, so the ESS estimates in \cref{fig:ess} suggest that the estimated moments of the posterior cannot be trusted, even with $2\times 10^5$ samples computed for this test case.

The bars in \cref{fig:est_expansion} indicate 3 times the length of the standard deviation. The standard deviation is estimated from the samples of the posterior, in both methods. We emphasize that the MCMC estimates of the standard deviation cannot be trusted due to the low ESS. Regarding the VED estimates of the standard deviation, we see that the confidence in the estimates reduces for higher modes, i.e., larger indices. This behavior is also reported in other works, e.g. see \cite{2111.15620,1930-8337_2016_4_1007}. In addition, since samples provided by the VED posterior predictive are independent samples, the estimated variance can be trusted, up to the accuracy of the VED approximation.

\begin{figure}
   \centering
\begin{tabular}{cccc}
\includegraphics[width=3cm]{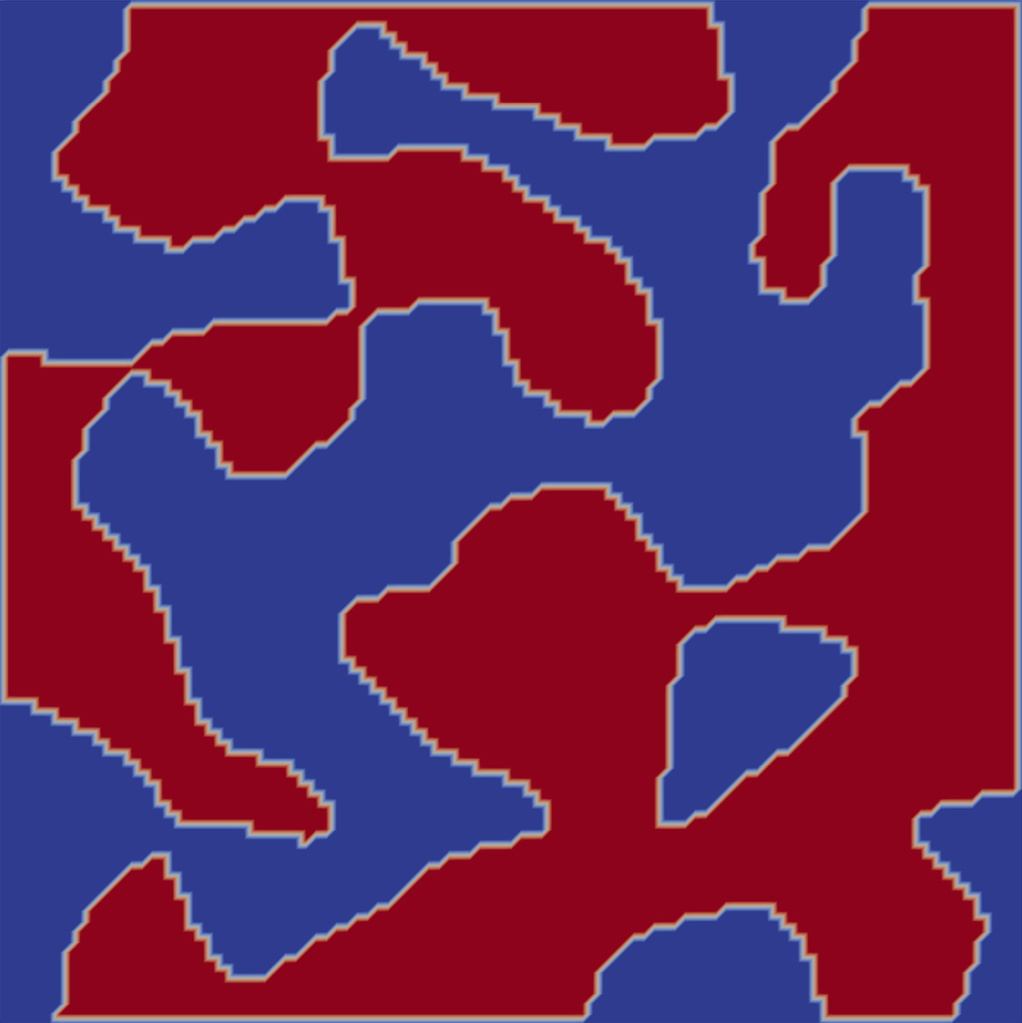}&
\includegraphics[width=3cm]{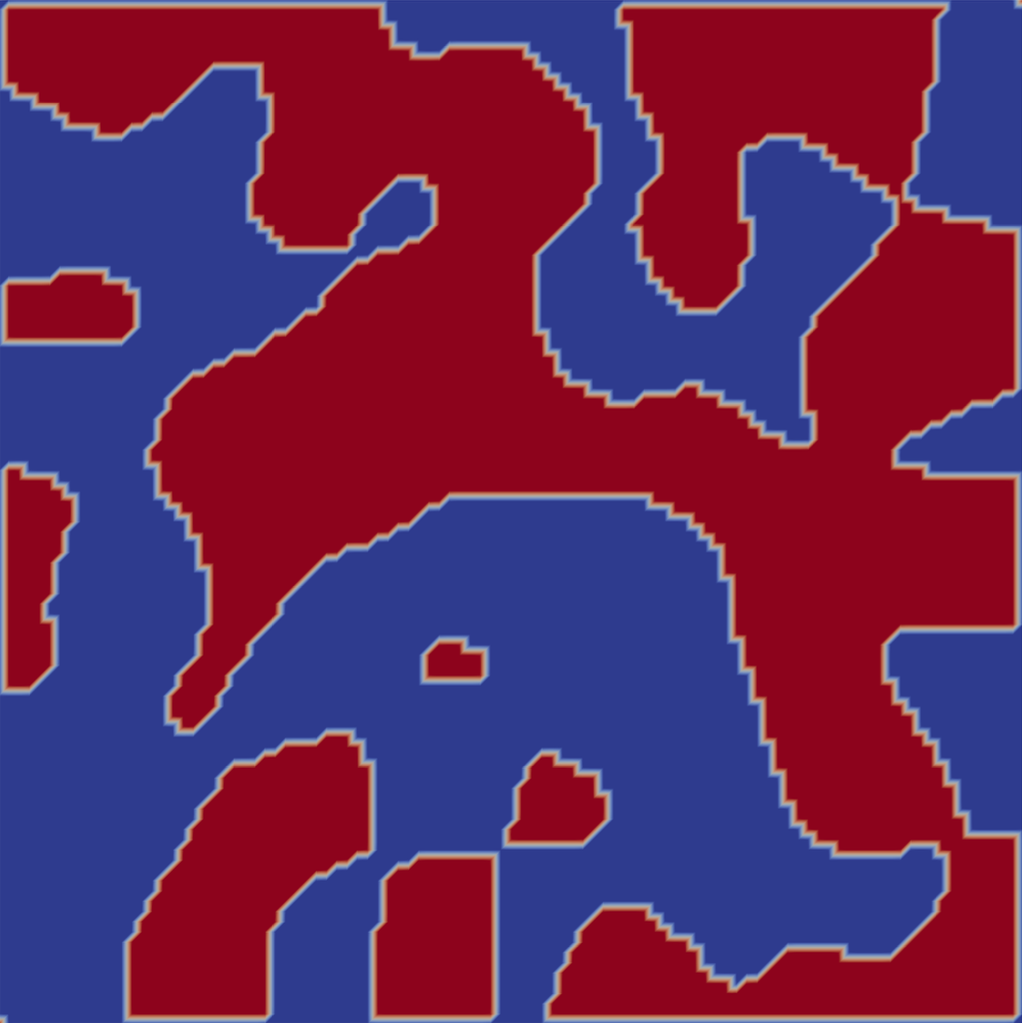}&
\includegraphics[width=3cm]{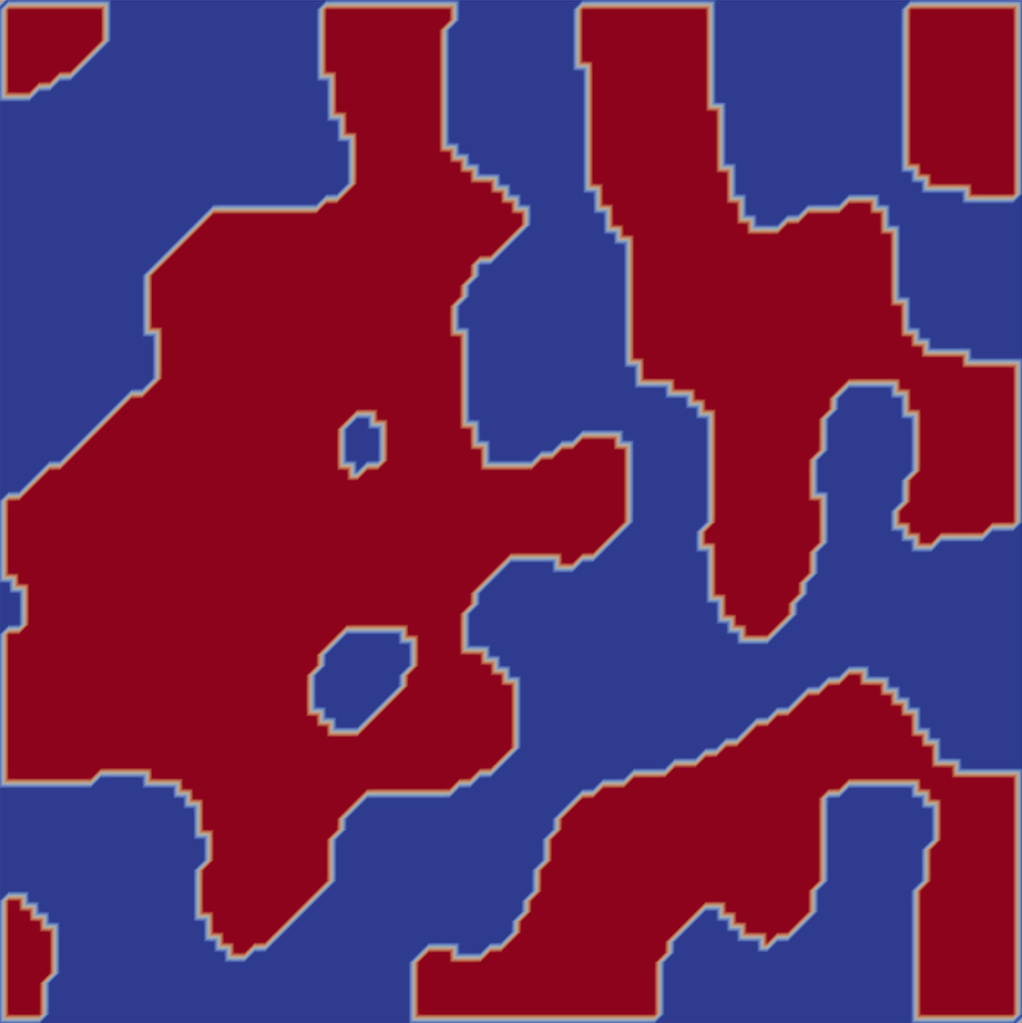}&
\includegraphics[width=3cm]{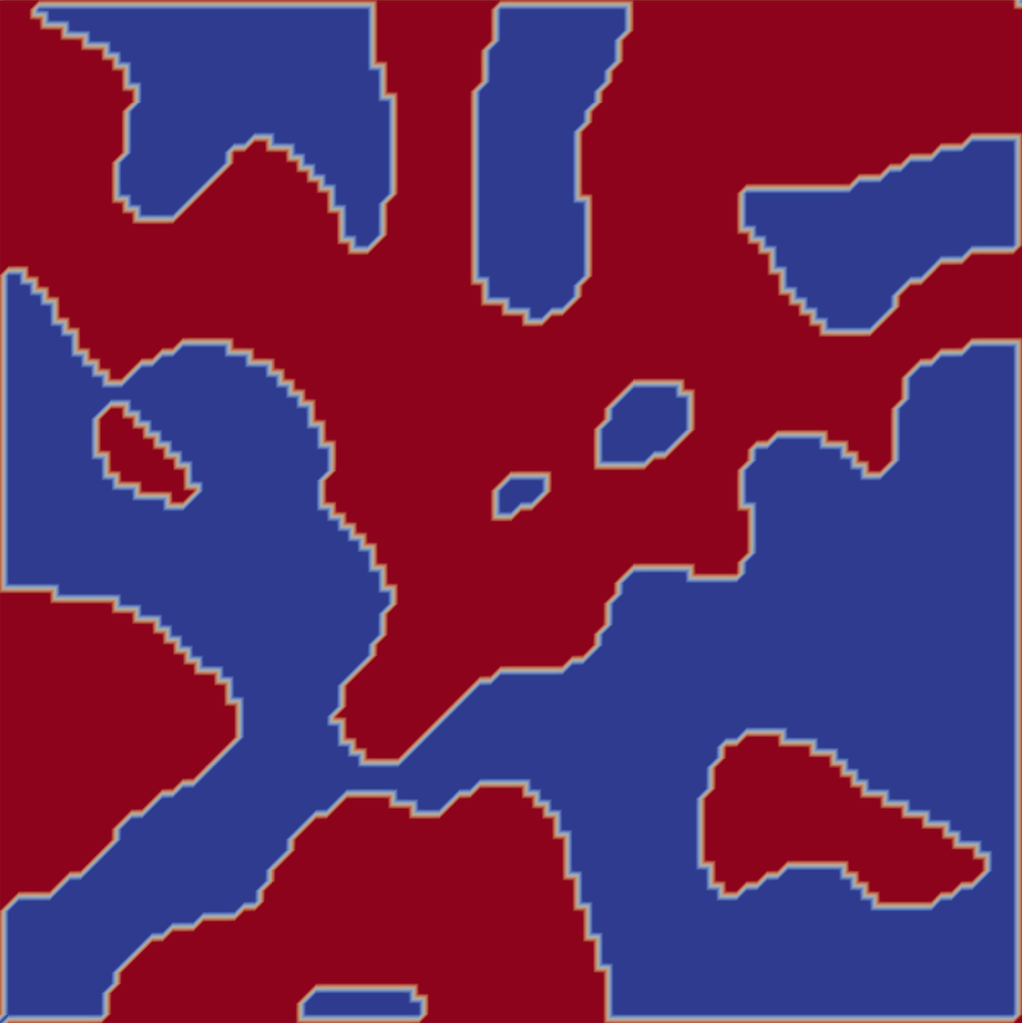} \\
\multicolumn{4}{c}{(a) samples from prior distribution
visualized on the conductivity field} \\
\includegraphics[width=3cm]{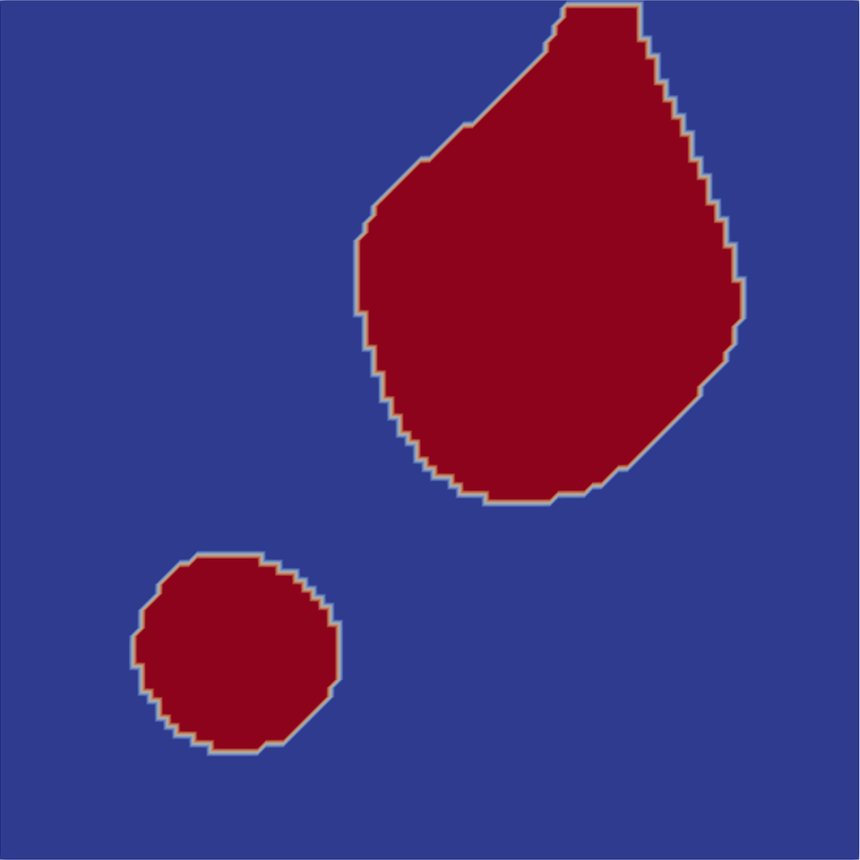}&
\includegraphics[width=3cm]{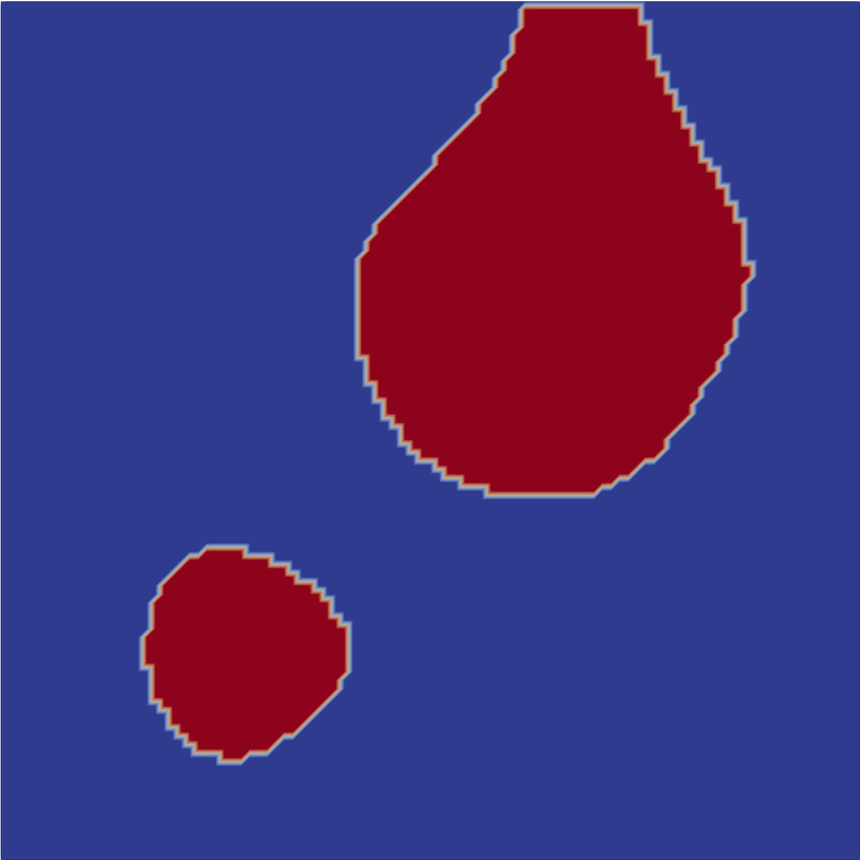}&
\includegraphics[width=3cm]{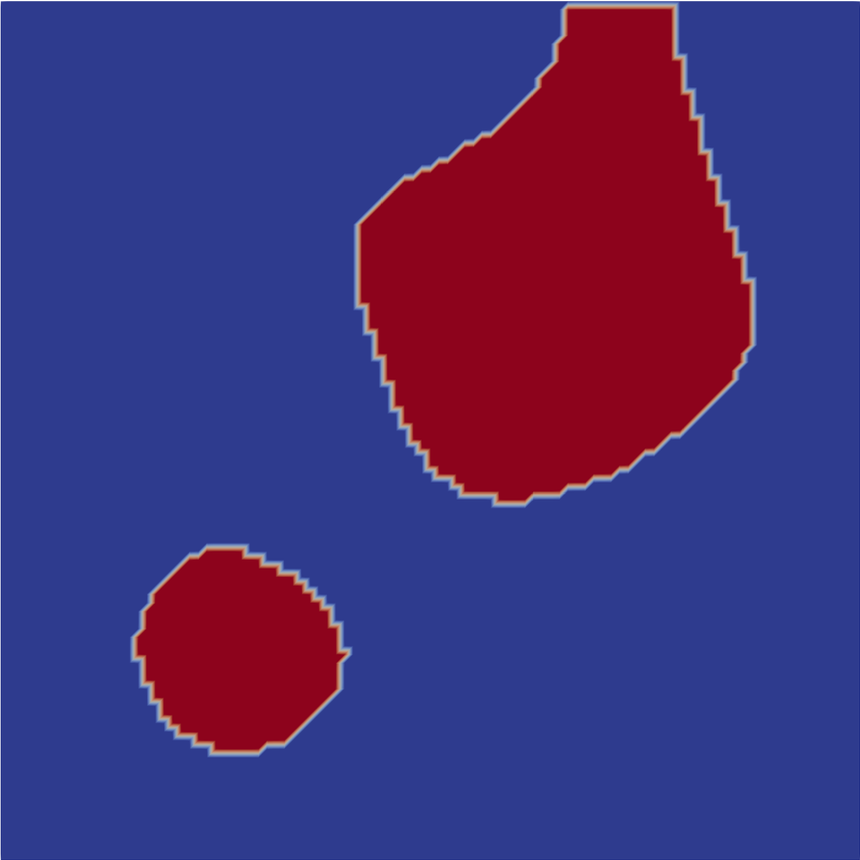}&
\includegraphics[width=3cm]{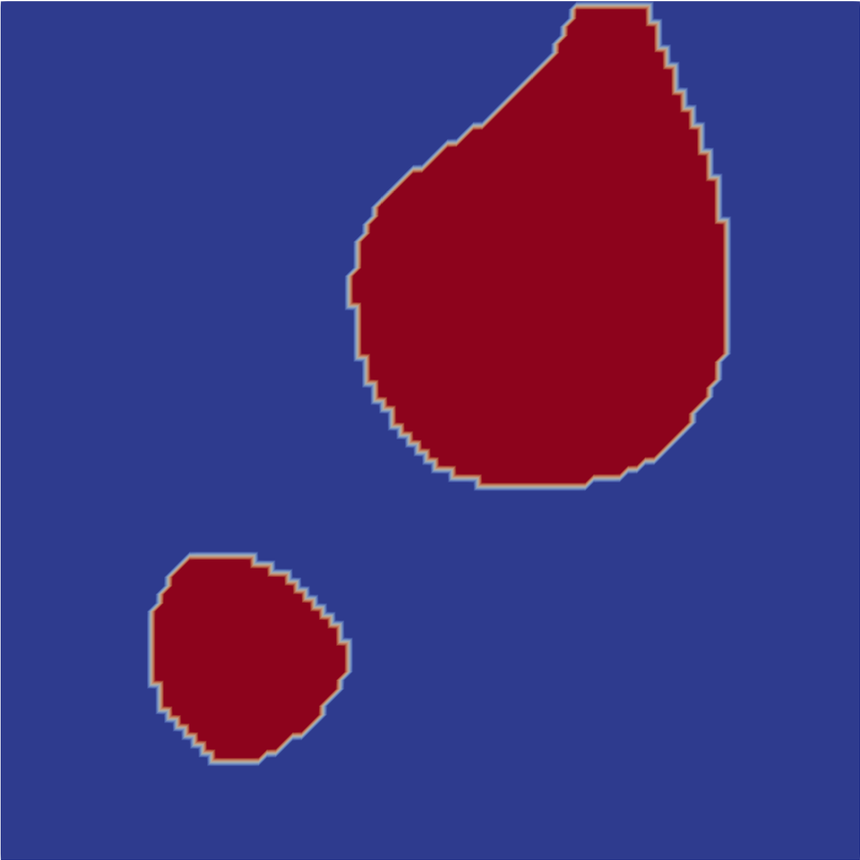} \\
\multicolumn{4}{c}{(b) samples from the VED posterior predictive visualized on the conductivity field}\\
\includegraphics[width=3cm]{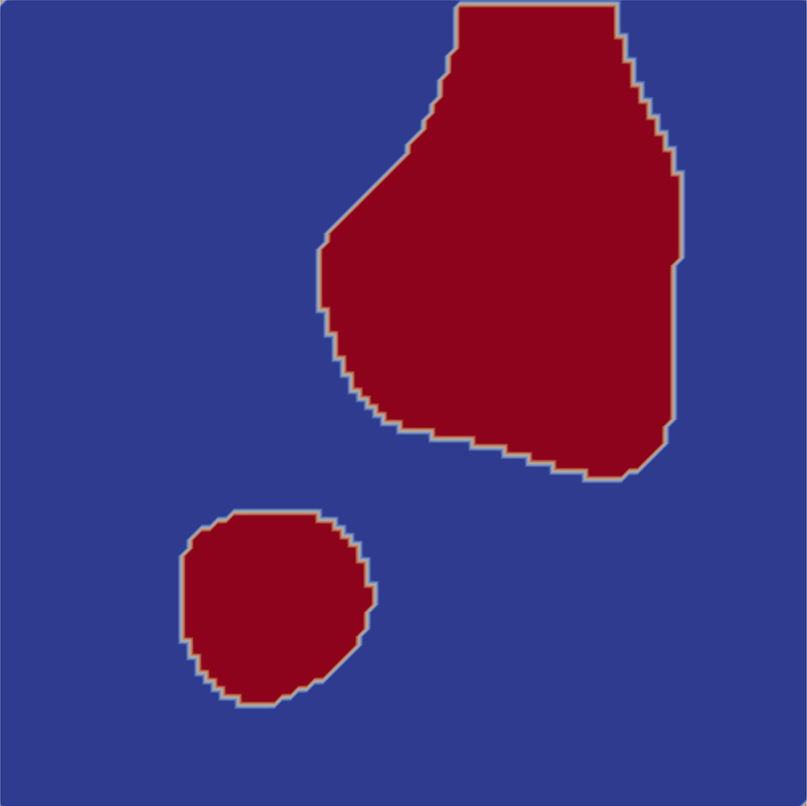}&
\includegraphics[width=3cm]{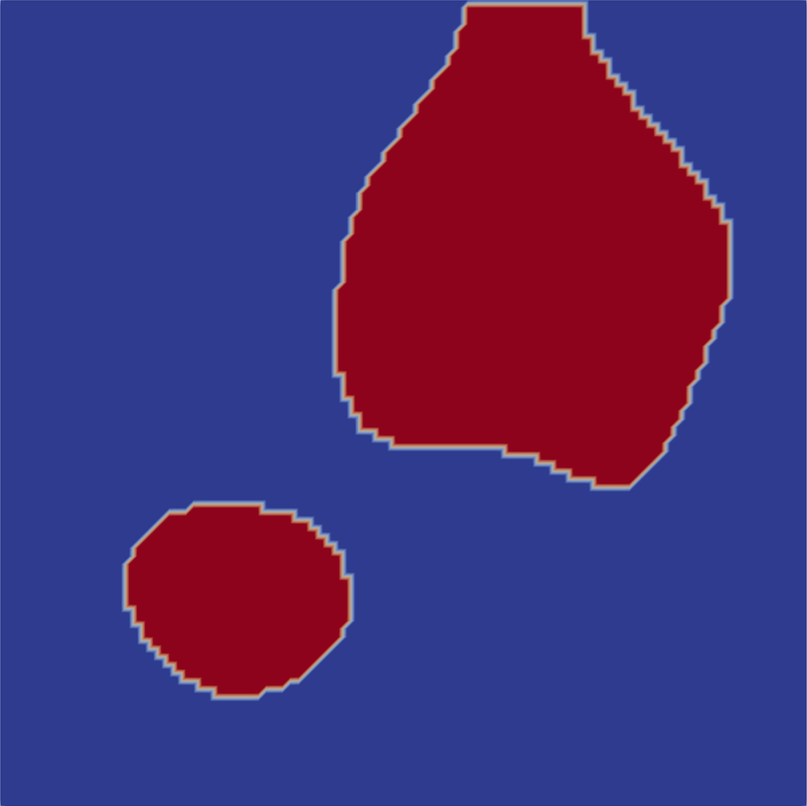}&
\includegraphics[width=3cm]{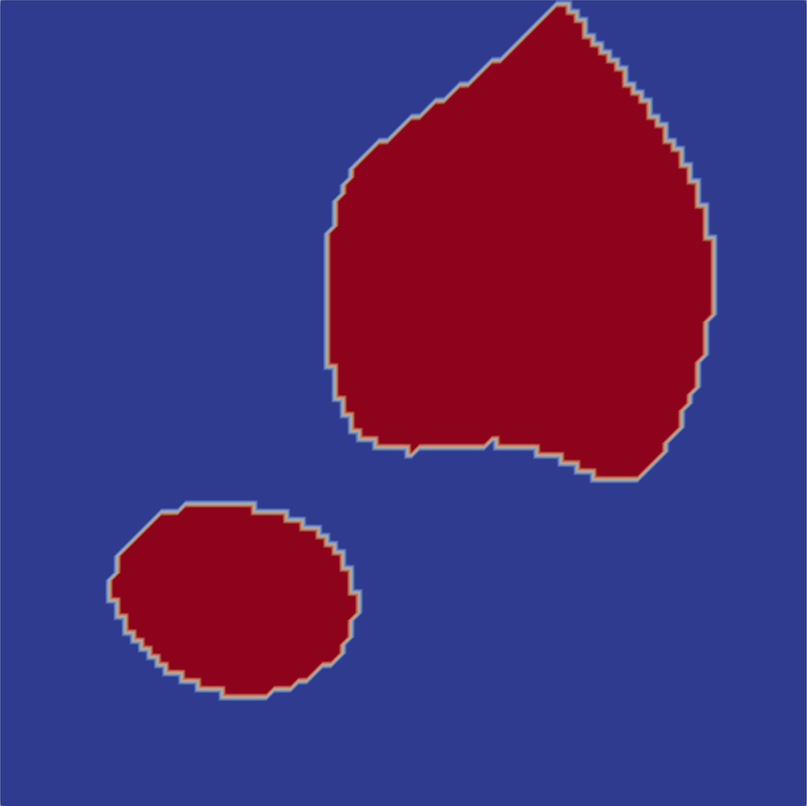}&
\includegraphics[width=3cm]{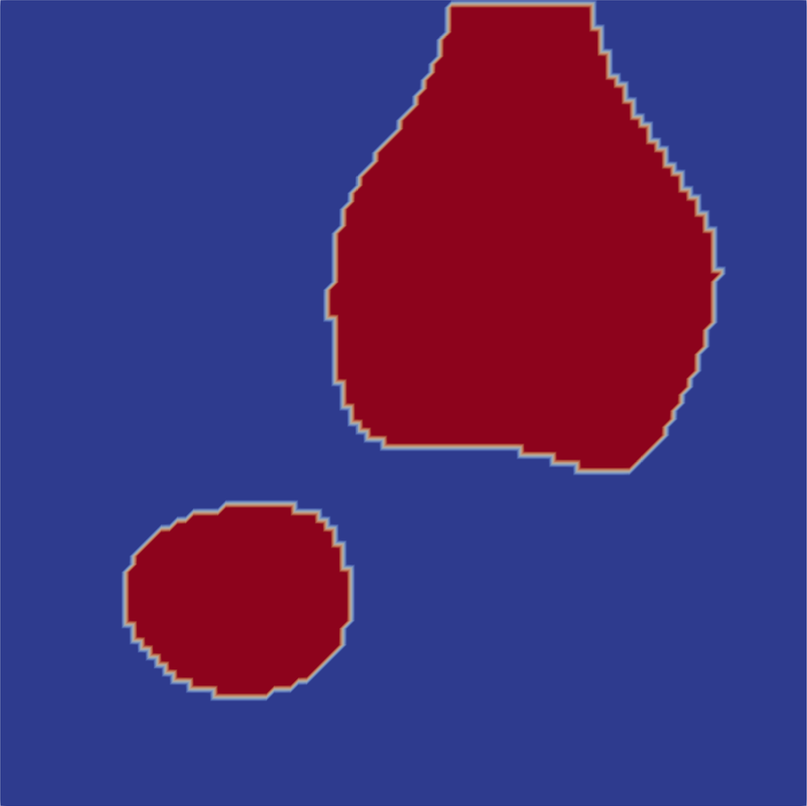}\\
\multicolumn{4}{c}{(c) MCMC samples from posterior predictive
visualized on the conductivity field}
\end{tabular}
\caption{(a) Samples from the prior that were obtained by \cref{eq:KL} and then \cref{eq:prior_conductivity} for random samples of the parameter vector $\bfx$. (b) Samples from the VED posterior predictive distribution mapped onto the conductivity field. (c) MCMC samples from the posterior predictive distribution mapped onto the conductivity field.}
\label{fig:prior-post}
\end{figure}

Next, we provide some samples for visualization.  In the first row of \cref{fig:prior-post}, we provide samples from the prior distribution that were obtained by mapping random samples of the parameter vector $\bfx^j$ onto the conductivity field $\bfy^j$ following \cref{eq:KL} and then \cref{eq:prior_conductivity}.  Then in the second and third rows of \cref{fig:prior-post}, we provide conductivity fields corresponding to samples from the VED posterior predictive distribution and MCMC samples from the posterior predictive distribution, respectively.
In terms of the conductivity field, we observe comparable samples.

\begin{figure}
   \centering
\begin{tabular}{cccc}
VED mean & VED variance & MCMC mean & MCMC variance \\
 \includegraphics[width=0.22\textwidth]{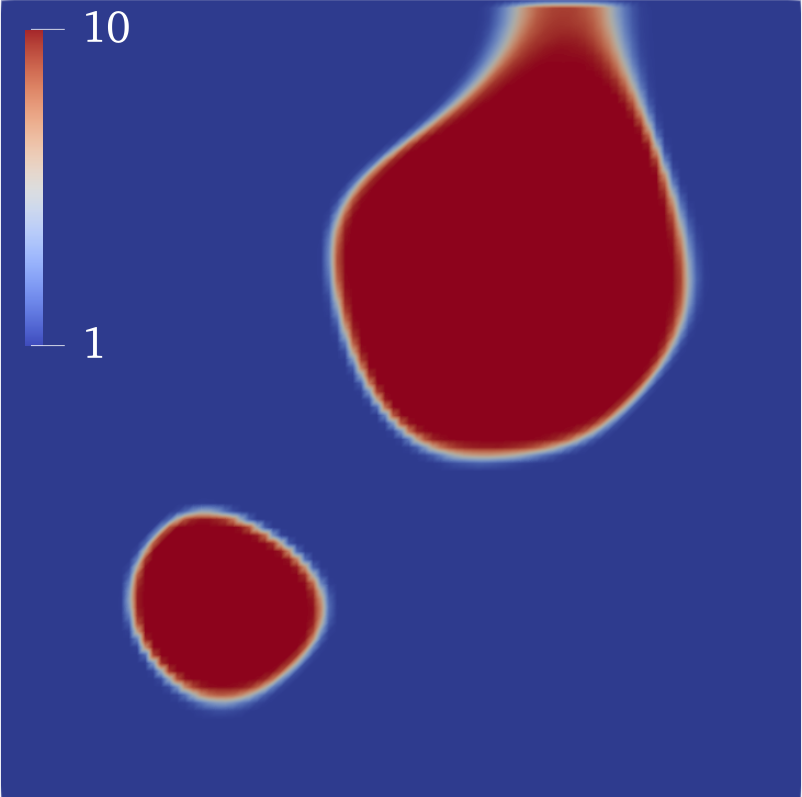}&
\includegraphics[width=0.22\textwidth]{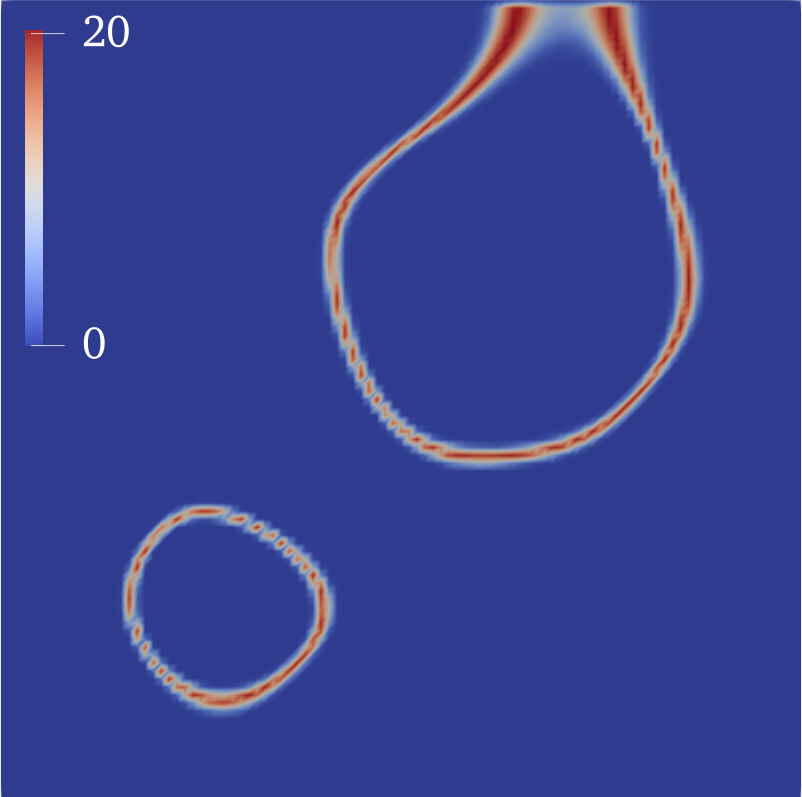} &
\includegraphics[width=0.22\textwidth]{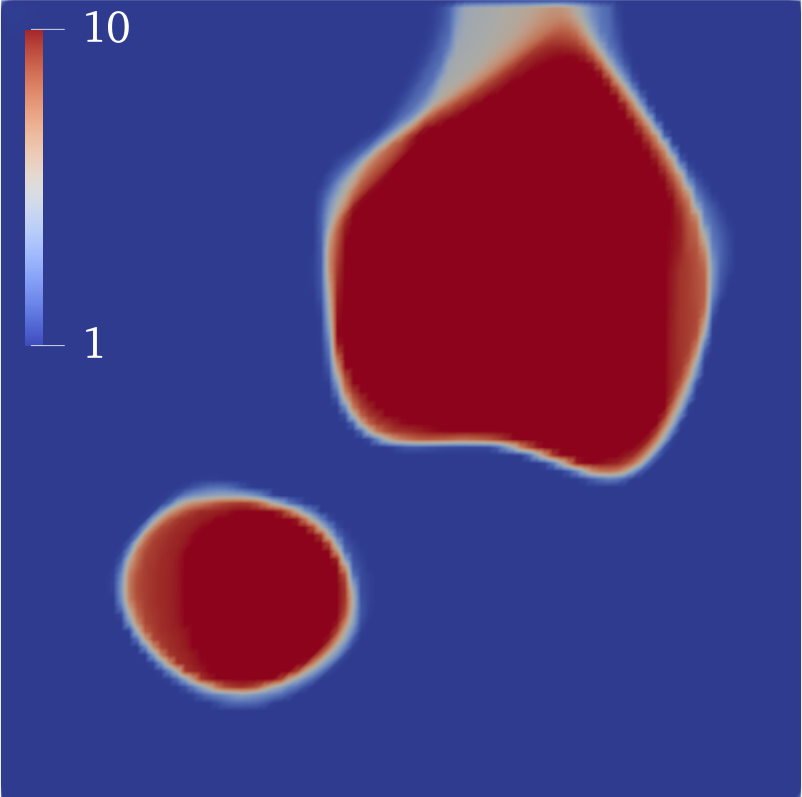} &
\includegraphics[width=0.22\textwidth]{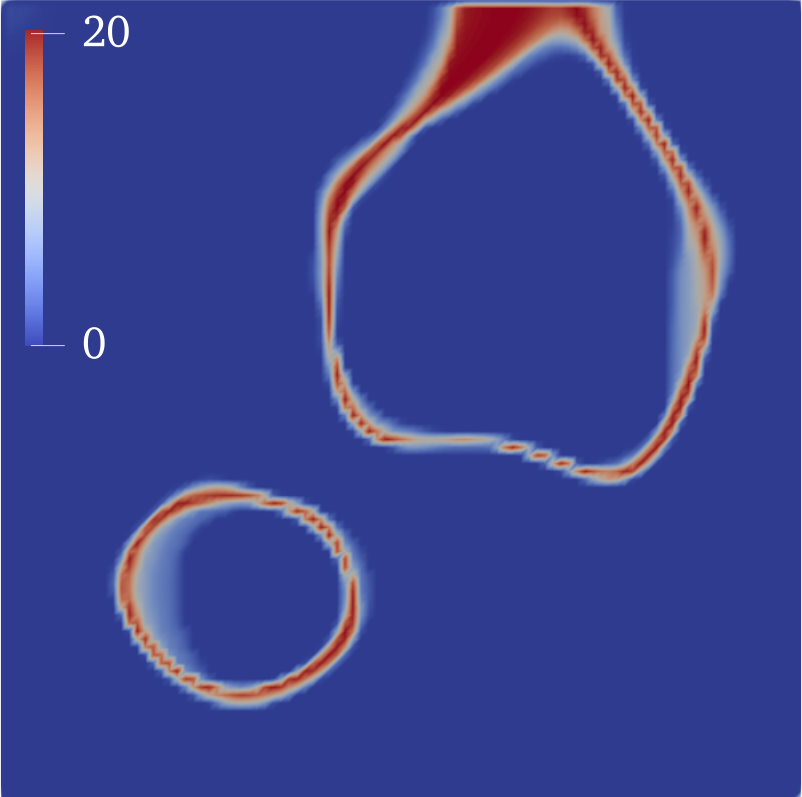}\\
\end{tabular}
\caption{Means and variances of the conductivity fields for samples from the VED posterior predictive (left) and MCMC samples from the posterior predictive (right). }
\label{fig:hydraulic_estimates}
\end{figure}

To visualize the uncertainty in the QoI, one approach is to compute the variance of the posterior-predictive for each expansion coefficient $X_i$ in \cref{eq:KL} (see \cref{fig:est_expansion}). Another approach is to compute the hydraulic conductivity for all samples from the poster predictive distribution using the transformation \cref{eq:prior_conductivity}, and then estimate means and variances in the conductivity space.  These results are provided in \cref{fig:hydraulic_estimates}.
Both methods correctly identify the two inclusions in the domain. The variance in both methods indicates a point-wise level of uncertainty in the estimation. Both methods show high certainty away from the boundaries of the inclusion and low certainty at the location of the boundaries. In addition, we see uncertainty at the boundary of the domain with the Newman boundary condition. However, as discussed above we cannot trust the uncertainty estimates for the MCMC method.

Finally, we compare the computational efficiency of the two methods by providing CPU times for computing $1000$ samples and $10$ independent samples from the posterior predictive in \cref{tab:time_compare}. All experiments were carried out on a MacBook Pro with an Apple M1 Pro processor. We notice that the ability to draw independent samples from the posterior-predictive provides significant computational gains. Recall that the number of samples required to train the neural network is $10^4$. This is orders of magnitude smaller than the number of samples required from the MCMC method to construct a representative set of independent samples from the posterior predictive distribution.

\begin{table}
\begin{center}
\begin{tabular}{ | m{7em} | m{5cm}| m{5cm} | }
 \hline
  & Elapsed time (s) for computing 1000 samples & Elapsed time (s) for computing 10 \textit{independent} samples \\
 \hline
 MCMC (pCN) & 270 & $2.5\times 10^5$ \\
 \hline
 VED sampling & 20 & 0.02 \\
 \hline
speed-up & 13.5 & $1.2\times 10^8$ \\
\hline
\end{tabular}
\end{center}
\caption{CPU times for sampling from the posterior predictive using MCM versus using VED sampling. The speed-up in terms of obtaining independent samples is significant.}  \label{tab:time_compare}
\end{table}

\section{Conclusions}
\label{sec:conclusions}
In this work, we describe VED networks for the direct estimation of QoI and their uncertainties for inverse problems, given observations. Following a supervised learning approach, we use training data to learn parameters defining a VED network, where an encoder maps an observation to a distribution over the latent space and a decoder maps samples from the latent space to distributions over the QoI output space. By leveraging probability distributions in both the latent space and in the posterior predictive, VED networks can be used to efficiently generate independent samples from the posterior predictive, thereby enabling efficient UQ for the QoI without requiring the use of expensive MCMC methods.  The potential benefits of our approach for goal-oriented UQ are demonstrated by the numerical experiments for two scenarios.  In \Cref{sub:experiment1}, we demonstrate the use of VED networks for UQ of the regularization parameter for total variation for X-ray CT reconstruction, and in \Cref{sub:experiment2}, we demonstrate the ability of trained VEDs to not only estimate locations but also uncertainties of discontinuities in nonlinear hydraulic tomography reconstructions.

\section*{Acknowledgement}
We would like to express our sincere gratitude to Dr. Heydar M. Afkham for his invaluable contribution in designing the architecture of neural networks, which has greatly enhanced the accuracy and efficiency of the models in this paper.

\printbibliography

@article{afkham2021learning,
  title={Learning regularization parameters of inverse problems via deep neural networks},
  author={Afkham, Babak Maboudi and Chung, Julianne and Chung, Matthias},
  journal={Inverse Problems},
  volume={37},
  number={10},
  pages={105017},
  year={2021},
  publisher={IOP Publishing}
}

@article{kingma2013auto,
  title={Auto-encoding variational Bayes},
  author={Kingma, Diederik P and Welling, Max},
  journal={arXiv preprint arXiv:1312.6114},
  year={2013}
}

@incollection{goh2019solving,
  title={Solving {B}ayesian inverse problems via variational autoencoders},
  author={Goh, Hwan and Sheriffdeen, Sheroze and Wittmer, Jonathan and Bui-Thanh, Tan},
  booktitle={Mathematical and Scientific Machine Learning},
  pages={386--425},
  year={2022},
  organization={PMLR}
}

@article{liu2021machine,
  title={Machine-learning-based prediction of regularization parameters for seismic inverse problems},
  author={Liu, Shihuan and Zhang, Jiashu},
  journal={Acta Geophysica},
  pages={1--12},
  year={2021},
  publisher={Springer}
}

@article{kingma2014adam,
  title={Adam: A method for stochastic optimization},
  author={Kingma, Diederik P and Ba, Jimmy},
  journal={arXiv preprint arXiv:1412.6980},
  year={2014}
}

@article{hammernik2018learning,
  title={Learning a variational network for reconstruction of accelerated MRI data},
  author={Hammernik, Kerstin and Klatzer, Teresa and Kobler, Erich and Recht, Michael P and Sodickson, Daniel K and Pock, Thomas and Knoll, Florian},
  journal={Magnetic Resonance in Medicine},
  volume={79},
  number={6},
  pages={3055--3071},
  year={2018},
  publisher={Wiley Online Library}
}

@article{lucas2018using,
  title={Using deep neural networks for inverse problems in imaging: beyond analytical methods},
  author={Lucas, Alice and Iliadis, Michael and Molina, Rafael and Katsaggelos, Aggelos K},
  journal={IEEE Signal Processing Magazine},
  volume={35},
  number={1},
  pages={20--36},
  year={2018},
  publisher={IEEE}
}

@book{calvetti2007introduction,
  title={An Introduction to Bayesian Scientific Computing: Ten Lectures on Subjective Computing},
  author={Calvetti, Daniela and Somersalo, Erkki},
  volume={2},
  year={2007},
  publisher={Springer Science \& Business Media},
%   address = {New York}
}

@article{farquharson2004comparison,
  title={A comparison of automatic techniques for estimating the regularization parameter in non-linear inverse problems},
  author={Farquharson, Colin G and Oldenburg, Douglas W},
  journal={Geophysical Journal International},
  volume={156},
  number={3},
  pages={411--425},
  year={2004},
  publisher={Blackwell Science Ltd Oxford, UK}
}

@article{vogel1996non,
  title={Non-convergence of the L-curve regularization parameter selection method},
  author={Vogel, Curtis R},
  journal={Inverse Problems},
  volume={12},
  number={4},
  pages={535},
  year={1996},
  publisher={IOP Publishing}
}

@article{mead2008newton,
  title={A Newton root-finding algorithm for estimating the regularization parameter for solving ill-conditioned least squares problems},
  author={Mead, Jodi L and Renaut, Rosemary A},
  journal={Inverse Problems},
  volume={25},
  number={2},
  pages={025002},
  year={2008},
  publisher={IOP Publishing}
}

@article{galatsanos1992methods,
  title={Methods for choosing the regularization parameter and estimating the noise variance in image restoration and their relation},
  author={Galatsanos, Nikolas P and Katsaggelos, Aggelos K},
  journal={IEEE Transactions on Image Processing},
  volume={1},
  number={3},
  pages={322--336},
  year={1992}
}

@book{hansen2010discrete,
  title={Discrete Inverse Problems: Insight and Algorithms},
  author={Hansen, Per Christian},
  year={2010},
  publisher={SIAM}
}

@online{randomSheppLogan,
  author = {Matthias Chung},
  title = {Random-Shepp-Logan-Phantom},
  year = 2020,
  url = {https://github.com/matthiaschung/Random-Shepp-Logan-Phantom},
  urldate = {2020-12-14}
}

@article{ruthotto2018optimal,
  title={Optimal experimental design for inverse problems with state constraints},
  author={Ruthotto, Lars and Chung, Julianne and Chung, Matthias},
  journal={SIAM Journal on Scientific Computing},
  volume={40},
  number={4},
  pages={B1080--B1100},
  year={2018},
  publisher={SIAM}
}

@article{haber2000gcv,
  title={A GCV based method for nonlinear ill-posed problems},
  author={Haber, Eldad and Oldenburg, Douglas},
  journal={Computational Geosciences},
  volume={4},
  number={1},
  pages={41--63},
  year={2000},
  publisher={Springer}
}

@article{de2016machine,
  title={A machine learning approach to optimal {T}ikhonov regularization I: affine manifolds},
  author={De Vito, Ernesto and Fornasier, Massimo and Naumova, Valeriya},
  journal={Analysis and Applications},
  volume={20},
  number={02},
  pages={353--400},
  year={2022},
  publisher={World Scientific}
}

@article{li2020nett,
  title={NETT: Solving inverse problems with deep neural networks},
  author={Li, Housen and Schwab, Johannes and Antholzer, Stephan and Haltmeier, Markus},
  journal={Inverse Problems},
  volume={36},
  number={6},
  pages={065005},
  year={2020},
  publisher={IOP Publishing}
}

@article{1930-8337_2016_4_1007,
title = {The Bayesian formulation of EIT: Analysis and algorithms},
journal = {Inverse Problems \& Imaging},
volume = {10},
number = {4},
pages = {1007},
year = {2016},
% issn = {1930-8337},
% doi = {10.3934/ipi.2016030},
% url = {http://aimsciences.org//article/id/d9e2dc77-0da0-4532-95df-b2c6e5a9eece},
author = {Matthew M.  Dunlop and Andrew M.  Stuart}
}

@inproceedings{bora2017compressed,
  title={Compressed sensing using generative models},
  author={Bora, Ashish and Jalal, Ajil and Price, Eric and Dimakis, Alexandros G},
  booktitle={International Conference on Machine Learning},
  pages={537--546},
  year={2017},
  organization={PMLR}
}

@article{peng2020solving,
  title={Solving inverse problems via auto-encoders},
  author={Peng, Pei and Jalali, Shirin and Yuan, Xin},
  journal={IEEE Journal on Selected Areas in Information Theory},
  volume={1},
  number={1},
  pages={312--323},
  year={2020},
  publisher={IEEE}
}

@article{hinton2006reducing,
  title={Reducing the dimensionality of data with neural networks},
  author={Hinton, Geoffrey E and Salakhutdinov, Ruslan R},
  journal={science},
  volume={313},
  number={5786},
  pages={504--507},
  year={2006},
  publisher={American Association for the Advancement of Science}
}

@inproceedings{salakhutdinov2007restricted,
  title={Restricted Boltzmann machines for collaborative filtering},
  author={Salakhutdinov, Ruslan and Mnih, Andriy and Hinton, Geoffrey},
  booktitle={Proceedings of the 24th international conference on Machine learning},
  pages={791--798},
  year={2007}
}

@inproceedings{torralba2008small,
  title={Small codes and large image databases for recognition},
  author={Torralba, Antonio and Fergus, Rob and Weiss, Yair},
  booktitle={2008 IEEE Conference on Computer Vision and Pattern Recognition},
  pages={1--8},
  year={2008},
  organization={IEEE}
}

@article{wang2016auto,
  title={Auto-encoder based dimensionality reduction},
  author={Wang, Yasi and Yao, Hongxun and Zhao, Sicheng},
  journal={Neurocomputing},
  volume={184},
  pages={232--242},
  year={2016},
  publisher={Elsevier}
}

@book{goodfellow2016deep,
  title={Deep learning},
  author={Goodfellow, Ian and Bengio, Yoshua and Courville, Aaron},
  year={2016},
  publisher={MIT press}
}

@article{goodfellow2020generative,
  title={Generative adversarial networks},
  author={Goodfellow, Ian and Pouget-Abadie, Jean and Mirza, Mehdi and Xu, Bing and Warde-Farley, David and Ozair, Sherjil and Courville, Aaron and Bengio, Yoshua},
  journal={Communications of the ACM},
  volume={63},
  number={11},
  pages={139--144},
  year={2020},
  publisher={ACM New York, NY, USA}
}

@article{cardiff2009bayesian,
  title={{B}ayesian inversion for facies detection: {A}n extensible level set framework},
  author={Cardiff, M and Kitanidis, PK},
  journal={Water Resources Research},
  volume={45},
  number={10},
  year={2009},
  publisher={Wiley Online Library}
}

@article{lee2013bayesian,
  title={{B}ayesian inversion with total variation prior for discrete geologic structure identification},
  author={Lee, J and Kitanidis, PK},
  journal={Water Resources Research},
  volume={49},
  number={11},
  pages={7658--7669},
  year={2013},
  publisher={Wiley Online Library}
}

@misc{2111.15620,
Author = {Reese, William and Saibaba, Arvind K. and Lee, Jonghyun},
Title = {{B}ayesian {L}evel {S}et {A}pproach for {I}nverse {P}roblems with {P}iecewise {C}onstant {R}econstructions},
Year = {2021},
Eprint = {arXiv:2111.15620},
}

@article{attia2018goal,
  title={Goal-oriented optimal design of experiments for large-scale {B}ayesian linear inverse problems},
  author={Attia, Ahmed and Alexanderian, Alen and Saibaba, Arvind K},
  journal={Inverse Problems},
  volume={34},
  number={9},
  pages={095009},
  year={2018},
  publisher={IOP Publishing}
}

@article{spantini2017goal,
  title={Goal-oriented optimal approximations of {B}ayesian linear inverse problems},
  author={Spantini, Alessio and Cui, Tiangang and Willcox, Karen and Tenorio, Luis and Marzouk, Youssef},
  journal={SIAM Journal on Scientific Computing},
  volume={39},
  number={5},
  pages={S167--S196},
  year={2017},
  publisher={SIAM}
}

@article{devore2019computing,
  title={Computing a quantity of interest from observational data},
  author={DeVore, Ronald and Foucart, Simon and Petrova, Guergana and Wojtaszczyk, Przemyslaw},
  journal={Constructive Approximation},
  volume={49},
  number={3},
  pages={461--508},
  year={2019},
  publisher={Springer}
}

@article{lieberman2013goal,
  title={Goal-oriented inference: {A}pproach, linear theory, and application to advection diffusion},
  author={Lieberman, Chad and Willcox, Karen},
  journal={SIAM REVIEW},
  volume={55},
  number={3},
  pages={493--519},
  year={2013},
  publisher={SIAM}
}

@article{lieberman2014nonlinear,
  title={Nonlinear goal-oriented {B}ayesian inference: application to carbon capture and storage},
  author={Lieberman, Chad and Willcox, Karen},
  journal={SIAM Journal on Scientific Computing},
  volume={36},
  number={3},
  pages={B427--B449},
  year={2014},
  publisher={SIAM}
}

@article{kullback1951information,
  title={On information and sufficiency},
  author={Kullback, Solomon and Leibler, Richard A},
  journal={The annals of mathematical statistics},
  volume={22},
  number={1},
  pages={79--86},
  year={1951},
  publisher={JSTOR}
}

@book{lemaire2013structural,
  title={Structural reliability},
  author={Lemaire, Maurice},
  year={2013},
  publisher={John Wiley \& Sons}
}

@article{10.1214/13-STS421,
author = {S. L. Cotter and G. O. Roberts and A. M. Stuart and D. White},
title = {{MCMC Methods for Functions: Modifying Old Algorithms to Make Them Faster}},
volume = {28},
journal = {Statistical Science},
number = {3},
publisher = {Institute of Mathematical Statistics},
pages = {424 -- 446},
year = {2013},
doi = {10.1214/13-STS421},
URL = {https://doi.org/10.1214/13-STS421}
}

@book{mcbook,
   author = {Art B. Owen},
   year = {2013},
   title = {{M}onte {C}arlo theory, methods and examples}
}

@article{dahl2017computing,
  title={Computing segmentations directly from x-ray projection data via parametric deformable curves},
  author={Dahl, Vedrana Andersen and Dahl, Anders Bjorholm and Hansen, Per Christian},
  journal={Measurement Science and Technology},
  volume={29},
  number={1},
  pages={014003},
  year={2017},
  publisher={IOP Publishing}
}

@article{afkham2023uncertainty,
  title={Uncertainty quantification of inclusion boundaries in the context of X-ray tomography},
  author={Afkham, Babak Maboudi and Dong, Yiqiu and Hansen, Per Christian},
  journal={SIAM/ASA Journal on Uncertainty Quantification},
  volume={11},
  number={1},
  pages={31--61},
  year={2023},
  publisher={SIAM}
}

@book{doi:10.1137/1.9781611976670,
author = {Hansen, Per Christian and J{\o}rgensen, Jakob and Lionheart, William R. B.},editor = {Per Christian Hansen and Jakob J{\o}rgensen and William R. B. Lionheart},
title = {{C}omputed {T}omography: {A}lgorithms, {I}nsight, and {J}ust {E}nough {T}heory},
publisher = {Society for Industrial and Applied Mathematics},
year = {2021},
doi = {10.1137/1.9781611976670},
address = {Philadelphia, PA},
URL = {https://epubs.siam.org/doi/abs/10.1137/1.9781611976670},
eprint = {https://epubs.siam.org/doi/pdf/10.1137/1.9781611976670}
}

@article{rudin1992nonlinear,
title = {{N}onlinear total variation based noise removal algorithms},
journal = {Physica D: Nonlinear Phenomena},
volume = {60},
number = {1},
pages = {259-268},
year = {1992},
issn = {0167-2789},
doi = {https://doi.org/10.1016/0167-2789(92)90242-F},
url = {https://www.sciencedirect.com/science/article/pii/016727899290242F},
author = {Leonid I. Rudin and Stanley Osher and Emad Fatemi}
}

@misc{https://doi.org/10.5281/zenodo.1254206,
  doi = {10.5281/ZENODO.1254206},
  url = {https://zenodo.org/record/1254206},
  author = {H\"{a}m\"{a}l\"{a}inen,  Keijo and Harhanen,  Lauri and Kallonen,  Aki and Kujanp\"{a}\"{a},  Antti and Niemi,  Esa and Siltanen,  Samuli},
  language = {en},
  title = {Tomographic X-Ray Data Of A Walnut},
  publisher = {Zenodo},
  year = {2015},
  copyright = {Creative Commons Attribution 4.0}
}

@article{hastie2015statistical,
  title={Statistical learning with sparsity},
  author={Hastie, Trevor and Tibshirani, Robert and Wainwright, Martin},
  journal={Monographs on statistics and applied probability},
  volume={143},
  pages={143},
  year={2015}
}

@article{chung2021efficient,
  title={Efficient learning methods for large-scale optimal inversion design},
  author={Chung, Julianne and Chung, Matthias and Gazzola, Silvia and Pasha, Mirjeta},
  journal = {Numerical Algebra, Control and Optimization},
  issn = {2155-3289},
  doi = {10.3934/naco.2022036},
  url = {/article/id/63930a2980e8bc15ce45a91f},
  year={2022}
}

@article{uribe2022hybrid,
  title={A hybrid Gibbs sampler for edge-preserving tomographic reconstruction with uncertain view angles},
  author={Uribe, Felipe and Bardsley, Johnathan M and Dong, Yiqiu and Hansen, Per Christian and Riis, Nicolai AB},
  journal={SIAM/ASA Journal on Uncertainty Quantification},
  volume={10},
  number={3},
  pages={1293--1320},
  year={2022},
  publisher={SIAM}
}
\end{document}